\documentclass[11pt]{article}
\usepackage{graphicx}
\def\L{{\bf L}}
\def\co{\overline{co}}

\def\pro{\diamondsuit}
\def\ve{\varepsilon}

\def\M{{\mathcal Y}}
\def\T{{\mathcal T}}
\def\I{{\mathcal I}}

\def\Q{{\mathcal Q}}
\def\C{{\mathcal C}}
\def\SS{s}
\def\RR{r}
\def\U{{\mathcal U}}

\def\g{{\bf g}}
\def\m{{\bf y}}
\def\u{{\bf u}}

\def\p{{\bf p}}
\def\y{{\wp}}
\def\v{{\bf v}}
\def\q{{\bf q}}
\def\w{{\bf w}}
\def\H{{\mathcal H}}
\def\F{{\mathcal F}}
\def\bV{{\bf V}}
\def\bP{{\bf P}}

\def\sqr#1#2{\vbox{\hrule height .#2pt
\hbox{\vrule width .#2pt height #1pt \kern #1pt
\vrule width .#2pt}\hrule height .#2pt }}
\def\square{\sqr74}
\def\endproof{\hphantom{MM}\hfill\llap{$\square$}\goodbreak}

\def\T{{\cal T}}
\def\O{{\cal O}}

\def\C{{\mathcal C}}
\def\R{I\!\!R}
\def\hh{I\!\!H}
\def\rr{I\!\!R}
\def\ll{I\!\!L}

\def\ov{\overline}

\def\forall{\hbox{for all }~}
\def\P{{\mathcal P}}

\def\vs{\vskip 2em}

\def\be{\begin{equation}}
\def\beq{\begin{equation}}
\def\bel{\begin{equation}\label}
\def\eeq{\end{equation}}

\newtheorem{thm}{Theorem}[section]

\newtheorem{cor}{Corollary}[section]
\newtheorem{lma}{Lemma}[section]
\newtheorem{prop}{Proposition}[section]
\newtheorem{remark}{Remark}[section]
\newtheorem{definition}{Definition}[section]

\begin{document}
\title{\bf Moving Constraints as Stabilizing Controls
in Classical Mechanics}
\vs
\author{Alberto Bressan$^{(*)}$ and Franco Rampazzo$^{(**)}$\\
(*)~Department of Mathematics, Penn State University, \\
University Park,
Pa.~16802, USA.\\
(**)  Dipartimento di Matematica Pura ed Applicata, Universit\`a di
Padova,\\
Padova  35141, Italy.\\
E-mails: bressan@math.psu.edu~ and  ~rampazzo@math.unipd.it\\
}
\maketitle

\begin{abstract}
The paper analyzes a Lagrangian system which is controlled by directly assigning
some of the coordinates as functions of time, by means of frictionless constraints.
In a natural system of coordinates, the equations of motions contain terms
which are linear or quadratic w.r.t.~time derivatives of the control functions.
After reviewing the basic equations, we explain the significance of
the quadratic terms, related to geodesics orthogonal to a given foliation.
We then study the problem of stabilization of the system to a given point,
by means of oscillating controls.  This problem is first reduced to the
weak stability for a related convex-valued differential inclusion,
then studied by
Lyapunov functions methods.
In the last sections, we illustrate the results by means
of various mechanical examples.
\end{abstract}


\section{Introduction}


A mechanical system can be controlled in
two fundamentally different ways.   In a commonly adopted
framework
\cite{BL, NS},
the controller modifies the time evolution of
the system by applying
additional  forces.     This  leads to a control problem in standard form,
where the time derivatives of the state variables depend continuously
on the control function.

In other situations, also physically realistic,
the controller acts on the system by directly assigning the
values of some of the coordinates,
by means of time dependent constraints.
 The evolution of the remaining coordinates can then
be determined by solving an ``impulsive" control system, where
the derivatives of the state variables depend (linearly or quadratically)
on the time derivative of the control
function.
This alternative point of view
was introduced, independently,  in \cite{AB1} and in \cite{Marle}.

Motivated by this second approach,
in  the present paper we study the
following problem  of Classical Mechanics:

 \noindent {\it Consider a system where the state\,\, space \,\,is a \,\,\,product
 $\Q\times\U$ of finite-dimensional manifolds
$\Q$ and $\U$. Assume that one can  prescribe
the motion $t\mapsto \u(t)\in\U$ of the second component, by means of
frictionless constraints.
 Given a point $(\bar \q,\bar \u)$,  can one stabilize the system
 at this point, by an oscillatory motion of the control $\u(\cdot)$
 around $\bar \u$ ?}

A well known example where stability is obtained by vibration
is provided by a pendulum whose suspension point can oscillate on a
vertical guide, as in Figure 1, left.  Calling $\theta$ the angle and $h$
the height of the pivot,
in this case we have $(\q,\u)=(\theta, h)\in S^1\times I$. Here
$S^1=[0,2\pi]$ with endpoints identified, and $I$ is an open interval.
If we take $\bar\q=\bar\theta =0$ as the (unstable)
upper vertical position of the pendulum,
it is well-known (see for example \cite{AKN, L1, L2} and references therein)
that this configuration can be made stable by rapidly oscillating the
pivot around a given value $\bar \u=\bar h$.
More generally, we will show that this system can be
asymptotically stabilized at any angle
$\bar\theta$ with $-\pi/2 < \bar\theta<\pi/2$, by a suitable
choice of the control function $t\mapsto h(t)=u(t)$.

\begin{figure}[h]
\centering
  \includegraphics[scale=0.50]{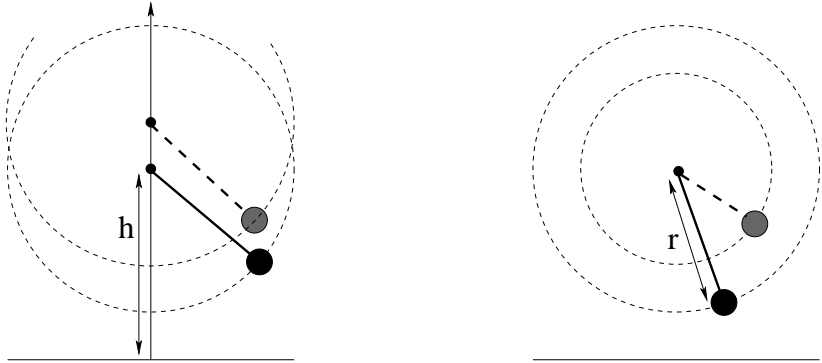}
\caption{Left: a pendulum with vertically moving pivot and fixed length. Right: a pendulum
with fixed pivot and variable length.} \end{figure}

On the other hand, consider the  variable length pendulum,
where the pivot is
fixed at the origin, but we can assign the radius of oscillation $r$
as function of time, see Figure 1, right.  The system is again
described by two coordinates
$(\q,\u)=(\theta, r)\in S^1\times I$.
However, in this case,
 the upright equilibrium position is {\em not} stabilizable by any
 oscillatory motion
of the radius $r(t)$  around a fixed value.

A major difference between these two systems is
that the equation of motion of the first one
 contains a quadratic term in the time derivative $\dot \u\doteq d\u/dt$.
On the other hand, the equation for the variable-length
pendulum is affine w.r.t.~the variable $\dot \u$.  Actually, the
explicit dependence on
 $\dot\u$ can be here entirely removed by a suitable change of coordinates.


 To understand the general problem, one has
 to consider two main issues.
  The former is {\bf geometric}, and  involves the
 {\em orthogonal curvature} of  the  foliation
 \bel{0.1}
 \Lambda\doteq \Big\{\Q\times\{\u\},\quad \u\in\U\Big\}.
 \eeq
Orthogonality is here defined  w.r.t.~the Riemannian
metric determined by the kinetic energy.
The orthogonal curvature
is a measure of how a geodesic,  which
is perpendicular to the leaf $\Q\times\{\u\}$ of the foliation at a given point
$(\q,\u)$,
 fails to remain perpendicular to the  other leaves it  meets.
   If this curvature is non-zero, then the dynamic equations
   for $\q$ and for the corresponding momentum $\p$ contain a quadratic term
in  the time derivative $\dot\u$ of the  control function.
This will be analyzed in detail in Part \ref{meccanica},
Sections \ref{ap1}, \ref{ap2}.

   The latter issue is {\bf analytical}, namely: how to
   exploit this curvature, i.e.~the quadratic terms in $\dot \u$, in in order
   to achieve stabilization.  This will be discussed in
   Part \ref{analisi} of this paper.  In particular, we study the set of
   solutions for a system
   with quadratic, unbounded, controls, making
   essential use of  reparametrization techniques.
   These, in turn, are combined with  arguments  involving
   Lyapunov functions for a convexified system.



The paper consists of three parts.
In order to keep our exposition as self-contained as possible,
in Part \ref{meccanica} we first describe
the mechanical model and derive the basic dynamical equations.
In Section \ref{mechanics} we recall
the classical equations of motion, in a Hamiltonian intrinsic
form, for time-dependent holonomic systems subject to non-conservative forces.
  Although this is a classical subject,
  which can be found in many text-books in Classical Mechanics, the purpose of
  this  first section  is to clarify concepts and notations
used in the remainder of the paper.
In Section~\ref{controlsec} we consider
a state space  $\M=\Q\times\U$ given by the product of two manifolds.
The  {\it controls}  will be  curves $t\mapsto \u(t)$ taking values in the
manifold $\U$. The main physical assumption we are making is that these controls
$\u(\cdot)$ are implemented by means of frictionless, time-dependent constraints.
One can then derive the equations of motion on the reduced state space
$\Q$, where the dynamics depends on
$\u$ and on its time derivative $\dot\u$, the latter dependence being
polynomial of degree two.
In Section \ref{loc}, we deduce the local expression of the control
equations
in a system of local coordinates adapted to the foliation
$\Lambda$ in (\ref{0.1}).
Section  \ref{ap1} contains a survey of some geometrical and functional
analytic results concerning the input-output map and the kinetic metric.
The main new result of Part~\ref{meccanica} appears in
Section \ref{ap2}, where  we present   a new
interpretation of the quadratic
dependence of the equations of motion on the derivative of the control
functions. Our characterization of the
quadratic coefficients is given in terms of the concatenation of two geodesics,
the second returning to the same leaf of the foliation
where the first one had started.
This generalizes to higher dimensions a result in
\cite{LR},  where the scalar control  case is considered. Finally,
Section~\ref{ap3} provides a variational characterization
of admissible control-trajectory pairs.

In Part \ref{analisi} we consider a
general nonlinear system where the right hand side is a quadratic
polynomial w.r.t.~the time derivatives of the control function.
\bel{0.2} \dot x= f(x)+\sum_{{\alpha}=1}^m g_{\alpha}(x)\,\dot u_{\alpha}
+\sum_{\alpha,\beta=1}^m
h_{{\alpha,\beta}}(x)\,\dot u_{\alpha}\dot u_{\beta}\,.
\eeq
Using  a re-parametrization technique, we show that the stabilization
problem for the impulsive control system (\ref{0.2}) can be reduced to proving a
weak stability property for a related differential inclusion with
compact, convex-valued right hand side:
\bel{0.3}\frac{d}{ds} x(s) \in F(x(s))\,,\eeq
$$
F(x)\doteq \co\left\{ f(x)\,w_0^2+
\sum_{{\alpha}=1}^m  g_{\alpha}(x) \,w_0 w_{\alpha}
+\sum_{\alpha,\beta=1}^m h_{{\alpha,\beta}}(x)\,w_{\alpha}w_{\beta}
~;~~~~w_0\in [0,1]\,,~~~\sum_{{\alpha}=0}^m
w_{\alpha}^2 = 1\right\}\,,$$
where $\overline{co}$ denotes a closed convex hull.
Theorems  \ref{5.1} and  \ref{thm5.2} relate the
weak (asymptotic) stabilizability of the differential inclusion
(\ref{0.3}) with the (asymptotic) stabilizability of the
impulsive control system (\ref{0.2}).

In practical cases, a direct analysis of the multifunction
$F$ may be difficult.  In Section \ref{sec:5},
in addition to (\ref{0.3})
we thus consider an auxiliary differential inclusion
of the form $\dot x\in G(x)$,
where the multifunction
$G$ is derived from
(\ref{0.2}) by neglecting all linear terms, i.e.~by formally setting
$g_\alpha\equiv 0$.
 We show that the weak  stability of this second
differential inclusion  still yields the relevant stabilization
properties for the original control system ({\ref{0.2}).
Motivated by \cite{Smi}, in Section \ref{sec:6} we
also show that the weak stability of the
differential inclusion can be established by looking at suitable selections.

In Part \ref{applicazioni meccaniche} we apply the previous analytic results
 to the problem of stabilization of mechanical systems,
 controlled by moving holonomic constraints.   Thanks to the particular
 structure of the
quadratic terms that appear in the equations of motion,
we show that in many cases one can construct a suitable
Lyapunov function, and thus
establish the desired stability properties.
The paper is then concluded with some examples, presented in Section \ref{examples}.

Throughout the paper, our focus is on systems in general form, where the equations of motion
depend quadratically on the time derivatives $\dot u_\alpha\,$.  In the special
case where the dependence is only linear, i.e.~$h_{\alpha,\beta}\equiv 0$ in
(\ref{0.2}), our results still apply; however,
controllability and stabilization are best studied by looking at Lie brackets
of the vector fields $f, g_\alpha$, using standard techniques of geometric control
theory  \cite{J, So2}.

In addition to \cite{AB1, Marle}, readers interested in the earlier developments
of the theory of
control of mechanical systems by moving constraints are referred
to \cite{AB2, CF, motta, Rampazzo1, Rampazzo2}.
A concise survey, also outlining
possible applications to swim-like motion in fluids,
 has recently appeared in \cite{B1}.  See also the lecture notes in
 \cite{RLect}.


\part{Time-dependent  holonomic constraints as controls}\label{meccanica}

\section{Review of the dynamical equations on the cotangent bundle}
\label{mechanics}\setcounter{equation}{0}

{\sc The Legendre-Fenchel transform}

Let $W$ be a finite-dimensional,
real vector space,
let $W^*$ be its dual space, and let $\langle\,\,,\,\rangle$ denote
the duality between $W$ and $W^*$. For every map $L: W\mapsto \rr$,
its Legendre-Fenchel transform
is defined as
\beq
\label{fenchel} L^*({\bf P}) \doteq \sup_{{\bV }\in W}
\{\langle {\bf P}, {\bV }\rangle - L({\bV })\}\, ,\eeq for
every ${\bf P}\in W^*$. $L^*$ is a convex map.  If $L$ is strictly convex, then $L^*$ is strictly convex, and $(L^*)^* = L$.

In the special case where $L$ is the sum of a positive
definite quadratic form and an affine function, the same is true of its transform
$L^*$. Moreover, the following facts are well-known:

\begin{itemize}
\item For every $\bP\in W^*$ there exists a unique $\bV  = i^L(\bP)\in W$
where the maximum on the right-hand side of (\ref{fenchel})
is  achieved. The map $ i^L:W^*\mapsto W$
is one-to-one,  affine, and satisfies
\beq
\label{inverse} (i^L)^{-1} = i^{L^*} .
\eeq
If $L$ is a quadratic form then $i^L$ and $i^{L^*}$ are linear.

\item By
identifying both $W$ and $W^*$ with $\rr^n$ (by the choice of a
basis on $W$ and of the dual basis on $W^*$) and using
$\frac{\partial L}{\partial \bV }$ and $\frac{\partial
L^*}{\partial \bP}$ to denote the gradients of $L$ and $L^*$,
respectively, one has $$ i^L(\bP) = \frac{\partial L^*}{\partial
\bP}(\bP)\qquad \qquad (i^{L})^{-1}(\bV )
=i^{L^*}(\bV ) = \frac{\partial L}{\partial
\bV }(\bV )\,.
$$
\end{itemize}

{\sc Holonomic mechanical systems }


Let $\M$ be a $d$-dimensional manifold, and call $T\M\doteq \{T_\m\,;~\m\in\M\}$,
~$T^*\M\doteq \{ T^*_\m\,;~\m\in\M\}$
its tangent and cotangent bundles.
If  $\mathcal W\subseteq\M$ is an open subset, and $Y:\mathcal W\mapsto \rr^d$ is a
coordinate
  chart, we recall that the corresponding {\it  bundle charts}
  $(Y,V)$ and $(Y,P)$ on $T\M$ and $T^*\M$,  are obtained
  by choosing the local  frames  $\left\{\frac{\partial}
  {\partial Y_{r}}\,;~~1\leq r\leq d
  \right\}$ and
    $\left\{dY_{r}\,;~~1\leq r\leq d\right\}$, respectively.

 By a {\em
 holonomic mechanical system} with $d$ degrees of freedom
 defined on a time interval $I\subset\rr$  we mean a pair
  $\Sigma = (\M,\T)$, where:
\begin{itemize}
\item $\M$ is a $d$-dimensional manifold; \item
$\T:\I\times T\M\mapsto \rr$ is a map, called the {\em kinetic energy
 of the system $\Sigma$},  such that
$ \T = \T_0 + \T_1 + \T_2\,, $, with
\beq\label{defT}\begin{array}{c}
 \T_0(t, {\m},{\bV }) \doteq \frac{1}{2} \g_0(t, \m)\,,
 \qquad\qquad
 \T_1(t, {\m},{\bV })\doteq \langle \g_1(t, \m) , {\bV }\rangle\,,
\\
\\
 \T_2(t, {\m},{\bV })\doteq \frac{1}{2}\g_2(t, \m)\big({\bV },{\bV }\big)\,.
\\
\end{array}
\eeq
\end{itemize}
In (\ref{defT}),  the $\g_{i}$ are smooth maps defined on $ I\times
\M$.  For each $(t,{\m})$, the map $\g_2(t,\m)$ is a positive-definite
quadratic form, $g_1(t,\m)$ is linear map, while
$g_0(t,\m)$ is a constant.

\vskip 1em

We say that the system $\Sigma=(\M,\T)$ is
{\em time-independent} if $\g_1(t,\m) = 0$, $\g_0
(t,\m) = 0$ for all $(t,\m)\in\I\times \M$, and $\g_2$ is independent of $t$.
Equivalently, $\T=\T_2$ and does not depend on  time.
In the frame-work of Lagrangian mechanics, this means
that the holonomic constraints determining the system $\Sigma$
are time-independent. In this case, the notion of mechanical
system coincides with that of Riemannian manifold, endowed with
the metric $\g_2$.

In terms  of a local bundle chart $(Y,V)$ on $T\M$, the maps in (\ref{defT})
take the form
\bel{kincoord}\T_0(t,Y,V) =
\frac{1}{2} {g_0},\qquad \T_1(t,Y,V)=
g_{r}V^{r},
\qquad \T_2(t,Y,V) =
\frac{1}{2}
{g}_{rs}V^{r}V^{s}.
\eeq
Here
the $d\times d$ matrix
$\big({g}_{r,s}\big)$,
the row vector $\big({{g}}_{r}\big)$, and the real number
$ g_0$ are the coordinate representations of $\g_2(t,{\m})$,
$\g_1(t,{\m})$,
 and  $ \g_0(t,{\m})$, respectively.  Here and in the sequel,
 it is understood that
 a summation should be performed over repeated indices.

Let us define the map $\T^*:I\times T^*\M\mapsto\rr$
as the Legendre transform of $\T$.
For every $(t,\m)\in I\times \M$, this means
 $$ \T^*(t,\m,{\bP}) = \sup_{{\bV }\in T_\m\M}\left\{
  \langle {\bP}, {\bV }\rangle -\T(t,\m,{\bV }) \right\}.$$
Similarly, we define  $\T^*_2$ as the Legendre
transform of $\T_2\,$.

The maps $\T^*$ and $\T_2^*$ will be also called the {\em
kinetic  Hamiltonians} corresponding to $\T$ and $\T_2$,
respectively. Accordingly, we shall use the notation
$$\H\doteq\T^*\qquad \H_2\doteq \T_2^*.$$

As in (\ref{inverse}),
for every $(t,\m)\in
  I\times \M$ we have the affine isomorphisms $i_{t,\m}^{\T}:
  T^*_{\m}\M \mapsto T_{\m}\M$ and
    $i_{t,\m}^{\T^*}:T_{\m}\M \mapsto T^*_{\m}\M$, defined as
 \beq \begin{array}{l}
 i_{t,\m}^{\T}\doteq i^{\T(t,\m,\cdot)}\,\,\qquad\hbox{and}\,\qquad\,
   i_{t,\m}^{\T^*} \doteq i^{\T^*(t,\m,\cdot)}\,.
  \end{array}
  \eeq
  respectively.
Entirely similar linear isomorphisms $i_{t,\m}^{\T_2}$ and $i_{t,\m}^{\T_2^*}$
can be defined in connection with the quadratic map $\T_2\,$.
According to (\ref{inverse}),   for all $(t,\m)\in
  I\times \M$ one has
  $$
  \left(i_{t,\m}^{\T} \right)^{-1} = i_{t,\m}^{\T^*} \quad \qquad
   \left(i_{t,\m}^{\T_2} \right)^{-1}  = i_{t,\m}^{\T_2^*}.
   $$

Using coordinates, if we use $ \big({g}^{r,s}\big)$ to denote
the inverse of the $d\times d$ matrix $ \big({g}_{r,s}\big)$
in (\ref{kincoord}),
then the affine isomorphisms
$i_{t,\m}^{\T}$
and $i_{t,\m}^{{\T}^*}=\left(i_{t,\m}^{\T}\right)^{-1}$
are given by
\beq\label{isoaff}
V^{s} = i_{t,\m}^{\T}(P_{1},\dots,P_d)\doteq
{{g}}^{r,s}\big(P_{r} -{{g}}_{r}\big)\qquad \qquad s=1,\dots,d\,,
\eeq
\beq
\label{isoaffinv} P_{s} =
i_{t,\m}^{\T^*}(V^{1},\dots,V^{d})\doteq
g_{r,s}V^{r} +{{g}}_{s}\qquad\qquad s=1,\dots,d
\,.
\eeq

\vskip 1em

By the identity
 $$
 {\mathcal H}(t, {\m}, {\bf P})~=~\T^*(t,\m,\bP)~=~
 \big\langle {\bf P}, i_{t,\m}^{\T}({\bf P})\big\rangle -
 \T\Big(t, {\m}, i_{t,\m}^{\T}({\bf P})\Big)\,,
 $$
 it straightforward to check that the Hamiltonian can be decomposed into
 a constant, a linear, and a quadratic part. Namely,
$ {\mathcal H} = {\mathcal H}_0 + {\mathcal H}_1 + {\mathcal H}_2\,$
 with
 $$ {\mathcal H}_2(t, {\m}, {\bf P}) ~=~
  \frac{1}{2}\,\g^2(t, {\m})({\bf P},{\bf P})~=~
  \frac{1}{2}\,\g_2(t,{\m})\Big( i_{t,\m}^{\T_2}({\bf P}),\,
  i_{t,\m}^{\T_2} ({\bf P})\Big)\,,
$$
 $$
 {\mathcal H}_1(t, {\m}, {\bf P}) ~=~
 \g^1(t, {\m})({\bf P})~=~ -\big\langle\g_1(t, {\m}),\,
  i_{t,\m}^{\T_2} ({\bf P})\big\rangle\,,
$$
$$
  {\mathcal H_0(t,\m)}~=~\frac{1}{2}\g^0(t, {\m})~ =~
  \frac{1}{2}\left[\g^2(t, {\m})\Big(\g_1(t, {\m}),\g_1(t, {\m})\Big)
   - \g_0(t, {\m})\right]\,.
  $$
We remark that, by our assumptions, the quadratic form
$\g^2(t,\m)$ is positive definite.
 \vskip 1em
Using local coordinates $(Y,P)$ the decomposition of the Hamiltonian function
takes the form
${\mathcal H}(t,Y,P) = {\mathcal H}_0 +
{\mathcal H_1} + {\mathcal H}_2\,$,
where
\bel{hamcoord}
{\mathcal H}_0\doteq\frac{1}{2}
{{g}}^{r,s}{{g}}_{r}{{g}}_{s} - \frac{g_0}{2}\,,\qquad
{\mathcal H}_1\doteq -
{{g}}^{r,s}P_{r}{{g}}_{s}\,,\qquad
{\mathcal H}_2\doteq \frac{1}{2}
{{g}}^{r,s}P_{r}P_{s}\,.
\eeq


{\sc Symplectic structure}

Consider again the $d$-dimensional manifold $\M$.  By
 $\omega_\M$ we shall denote  the canonical symplectic form
 on $T^*{\M}$. This
 is
 the (nondegenerate) $2$-form which, in any chart of $\M$ and the
corresponding bundle  chart of
 $T^*\mathcal \M$, is represented
   by the constant $2d\times 2d$ symplectic matrix
  $$
  S_d \doteq  \left(\begin{array}{ll} \,\,\,\,\,0_d&{\bf I}_d
  \\\,&\,\\-{\bf I}_d&0_d
 \end{array} \right).$$

For each $(\m,\bP)\in T^*\M$, by means of the
  canonical symplectic form $\omega_\M$ one can
  establish a linear  isomorphism
 $
 s_{\m,{\bf P}}:T_{\m,{\bf P}}^*\big( T^*\M\big)\mapsto
T_{\m,{\bf P}}\big( T^*\M\big).
$
This is uniquely defined
by setting
\bel{symplectic}
\langle \beta, {\bf w} \rangle = \omega_{\M}(\m,{\bf P})\Big(s_{\m,{\bf P}}(\beta), {\bf w}\Big) \,.
\eeq
 With reference to coordinates $(Y,P)$  ,
 if $\beta =
 \left({\beta_r} d Y^{r}+
\gamma_r dP_r\right)$,  one has
 $$
 s_{\m,{\bf P}}(\beta) =
 \gamma_r \frac{\partial }{\partial Y^{r}} -
 \beta_r\frac{\partial }{\partial P_{r}}\,,
 $$

For every function $\phi:I\times{T^*\M}\mapsto \rr$, differentiable
w.r.t.~the second variable for every $t\in I$,
 one can define a time-dependent vector field $X_{\phi}$ on
${T^*\M}$, called the\,\,\, {\em Hamiltonian vector field} associated to
$\phi$. For each $(t,{\m},{\bf P})\in I\times T^*\M$, this is defined by
$$
X_{\phi}(t, \m,{\bf P})\doteq s_{\m,{\bf P}}\big( d\phi(t,{\m},{\bf
P})\big)
$$
where $d\phi(t,\cdot)$ denotes the
differential of the map $\phi(t,\cdot):T^*\M\mapsto\rr $. In local coordinates,
one has
 $$
 X_{\phi} =
 \frac{\partial \phi}{\partial P_{r}}
 \frac{\partial }{\partial Y_{r}} -
 \frac{\partial \phi}{\partial Y^{r}}
 \frac{\partial }{\partial P_{r}}\,.
 $$

\vskip 2em

{\sc Forces}

Forces are represented by  vertical vector
fields on the cotangent bundle $T^*\M$.
We recall that, for $(\m,\bP)\in T^*\M$,
a vector $X\in T_{(\m,\bP)}(T^{*}\M)$ is called {\em vertical} if
$$
d{\Pi}(\m,\bP) \big(X\big) = 0.
$$
Here
${\Pi}: T^*\M \mapsto \M $
is the canonical projection of $T^*\M$ on $\M$ and
 $ d{\Pi}(\m,\bP)$
denotes its differential at $(\m,\bP)$. The subspace
  of vertical vectors   at $(\m,\bP)$ is thus the kernel of the map
  $ d{\Pi}(\m,\bP)$.
It will be called
the {\bf vertical tangent space} of $T^*\M$ at $(\m,\bP)$      and
denoted by $V_{(\m,\bP)}(T^{*}\M)$. The corresponding
fiber sub-bundle of $T(T^*\M)$
 is called the {\em vertical tangent bundle}
  of $T^*\M$ and  denoted by $V(T^*\M)$.

A vector field  ${\bf F}$ on $(T^*\M)$ is said to be
 {\em vertical }   if
${\bf F}(\m,\bP)\in V_{(\m,\bP)}(T^{*}\M)$ for each  $(\m,\bP)\in T^*\M$.
Using a canonical bundle chart of $T^*\M$, a vertical vector field is thus
 represented by
 a $2d$-dimensional column vector of the
form $(0,\ldots,0, F_1,\ldots, F_d)^\dag$, where the superscript~$^\dag$ denotes
transposition.

More generally, let $I$ be a real
interval. A function ${\bf F}:I\times {T^*\M}\mapsto V({T^*\M})$ such that, for
each $t\in I$, the map ${\bf F}(t,\cdot)$ is a vertical vector field on
$T^*\M$ will be called a {\em time-dependent vertical vector field}.
In order to retain its physical meaning,
${\bf F}$ will also be called a (possibly time-dependent)
{\em force} acting on $\M$ during the time-interval $I$.


  A time-independent force ${\bf F}$ is called {\em positional} if
${\bf F}(\m,\bP)= {\bf F(\m)}$, i.e.~if its values depend only on $\m$
and not on the value $\bP\in T_{\m}\,$.
Furthermore, a positional force ${\bf F}$ is called {\em conservative} if
there exists a  potential function $U:\M\mapsto \rr$
such that ${\bf F}$ is the Hamiltonian vector field corresponding
to $U$, namely
\beq\label{conservative} {\bf F} = X_{U}\,.
\eeq


Recalling the symplectic form $\omega_\M$
in (\ref{symplectic}), we now
introduce a formal notion of {\em power} of the force ${\bf F}$,
which will play an important role in Section \ref{controlsec}.
Let $(\m,\bP)\in T^*\M$,
   ${\bf F}\in V_{(\m,\bP)}(T^*\M)$, and  $\bV\in T(T^*\M)$. The quantity
   \bel{power}
   \omega_\M({\bf F}\,,\,\bV)
\eeq
   will be called the {\em power of ${\bf F}$ with respect to $\bV$}.
If $(Y, P)$ are canonical local coordinates
on $T^*\M$ and $${\bf F} = F_{r}\frac{\partial}{\partial P_{r}}\qquad\bV
   = V_{r}\frac{\partial}{\partial Y^{r}} +
    W_{r}\frac{\partial}{\partial P_{r}},
   $$
then one obtains the familiar expression for the power:
   $$
   \omega_\M({\bf F}\,,\,\bV) = F_{r}V_{r}\,.
   $$

\vskip1truecm

 {\sc The equation of motion}

Let $\Sigma=(\M,\T) $ be a mechanical system and let $\H$ be the corresponding
kinetic Hamiltonian. Let ${\bf F}$ be a force acting on $\M$. Then the
{\bf equation of motion}
 for the  mechanical system $\Sigma$
   subject to the force ${\bf F}$ is the
differential equation
\beq\label{motion}
 \frac{d}{dt}\left(\begin{array}{c}\m\\\bP\end{array}\right)= X_{\mathcal H} + {\bf F}\qquad t\in I
 ,\quad (\m,\bP)\in T^*\M\,,\eeq
where $X_{\mathcal H}$ is the Hamiltonian vector field
 associated to $\H$.

If  ${\bf F}$ is
a conservative force with potential function $U:\M\mapsto\R$,
we can consider the
standard Hamiltonian
 $$ H \doteq (\T-U)^*\,,
$$
defined as the Legendre transform of $\T-U$.
Then
(\ref{motion}) reduces to the usual Hamiltonian form
\beq\label{motioncons}
 \frac{d}{dt}\left(\begin{array}{c}\m\\\bP\end{array}\right)
 = X_H\qquad  t\in I\,.
 \eeq
 Indeed, one has $X_{\H} + X_{U} = X_{\H+U}$ and $\H+U =
(\T-U)^*=H\,$.

  With the usual notational
   conventions, the equation of motion (\ref{motion})
    in a local bundle chart takes the form

\beq\label{motioncoord}
  \left\{\begin{array}{l} \dot Y^r =
  \frac{\partial\H}{\partial P_{r}}\\
\,\qquad\qquad\qquad\qquad\qquad \qquad\qquad r = 1,\dots,d\,.\\
\dot P_{r} = - \frac{\partial\H}{\partial Y^{r}}
 + F_{r}  \end{array}\right.
 \eeq

Using the expressions (\ref{hamcoord}) for the Hamiltonian,
from (\ref{motioncoord}) we obtain
  \beq\label{motioncoordexpl}
  \left\{\begin{array}{l} \dot Y^{r} =
  {{g}}^{r,s}\Big(P_{s} -
 {{g}}_{s}\Big)\,,\\
\,\qquad
\\
\dot P_{r} =
-\frac{\partial g^{\ell s}}{\partial Y^r}
\left(\frac{1}{2} g_\ell g_s- P_\ell g_s +\frac{1}{2}
P_\ell P_s\right)
+ g^{\ell s}\,\frac{\partial g_s}{\partial Y^r} \left( P_\ell - g_\ell\right)
+ \frac{\partial g_0}{\partial Y^r} +
F_{r}\,.
 \end{array}
 \right.
 \eeq
 In the case of a time-independent system,
 these equations reduce to
  \beq\label{autonomous}
  \left\{\begin{array}{l} \dot Y^r =
  g^{r,s}P_{s}\,, \\
\,\qquad\\
\dot P_{r} = -\frac{1}{2}
\frac{\partial
g^{\ell s}}{\partial Y^{r}}P_\ell P_s + F_r\,.
\end{array}
 \right.
 \eeq
In particular, if the force is conservative,
(\ref{autonomous}) takes the familiar Hamiltonian form
 \beq\label{cons}
  \left\{\begin{array}{l} \dot Y^{r} = \frac{\partial H }{\partial P_{r}}
  \,\,\,\,\,\,\,\,\,\,\,\,\big(=
  g^{r,s}P_{s}\big)\,,\\
\,\qquad
\\
\dot P_{r} = -\frac{\partial H }{\partial Y^r}
\,\,\,\,\,\,\,\,\,\big(=-\frac{1}{2}
\frac{\partial g^{\ell s}}{\partial Y^r}P_{\ell}P_{s} -
\frac{\partial U}{\partial Y^{r}} \big)\,.
\end{array}
 \right.
 \eeq

\section{Time-dependent constraints as controls}
\label{controlsec}\setcounter{equation}{0}
 In this section we shall set up the general framework to
treat
 the situation where additional time-dependent holonomic constraints
 are regarded as {\it controls}.

 \vskip0.2truecm \subsection{ Structural assumptions} We shall consider a
mechanical system $\Sigma= (\M,\T)$
verifying the following assumptions:
 \begin{itemize}
 \item[{\bf 1)}]({\sc Product structure}). The state manifold $\M$ is a product
manifold, namely
\bel{prodmanif}\M = \Q\times\U\,,
\eeq
Here $\Q$ and $\U$, called the
{\it  reduced state space} and the {\it control space},
  are manifolds of dimension $N$ and $M$, respectively.

 \item[{\bf 2)}]({\sc Stationarity of the metric})\footnote{This
 assumption is made only for the sake of simplification, since
 it avoids  a double time-dependence: the structural
 one and the one due to the implementation of controls. Actually,
the situation where $\g$ has a general form can be treated as
well without significative additional difficulties.}
  The Lagrangian system $\Sigma$ is time-independent. Namely, the
  kinetic energy $\T$ is defined by a triple
   $\g=(\g_2,0,0)$ with
   $\g_2$ independent of time.


\item[{\bf 3)}]({\sc Regularity of the force}). The
external force ${\bf F}={\bf F}(t,\q,\u,\bP,\y)$ is a function measurable w.r.t.~$t$
and locally Lipschitz w.r.t.~all other variables.
\end{itemize}

By (\ref{prodmanif}), one has the natural identifications of
   $T^*(\Q\times\U)$,\,\, $T(T^*(\Q\times\U))$, \,\, and \,\, $V(T^*(\Q\times\U))$
with \,\, the products \,\, $T^*(\Q)\times T^*(\U)$,\,\, $T(T^*(\Q))\times
T(T^*(\U))$, \,\,and $V(T^*(\Q))\times V(T^*(\U))$, respectively.

By the second assumption, $\T = \T_2$, and,  for every $(t,\m)\in I\times\M$,
 the affine isomorphism  $i^{\T}_{t,\m}= i^{\T}_{\m}$ is a linear isomorphism,
 independent of time.
 In turn the Hamiltonian
is time-independent and its global and local expressions are given
by
 $$
  {\mathcal H}({\m},{\bf P}) =
  \frac{1}{2}\g^{-1}({\m})({\bf P},{\bf P})
  $$
 and
 $$
{\mathcal H}(Y,P) =\frac{1}{2}\sum_{r,s=1}^{N+M} {{
g}}^{r,s}(Y)P_{r}P_{s} \,\,,
$$
 respectively.
Moreover, for each ${\m}\in \M$, as soon as the vector spaces
$T_{\m}\M$ and $T_{\m}^*\M$ are endowed with the scalar products
defined by $\g({\m})$ and $\g^{-1}({\m})$, respectively, the
isomorphism ${i_{\m}^{\T}}$ is an isometry. In particular, one has
$$
\T\Big({\m},{i_{\m}^{\T}}({\bf P})\Big) = {\mathcal H}({\m},{\bf P})
$$
for all ${\m}\in \M$ and ${\bf P}\in T^*_{\m}\M$.

We observe that,
for every $(\q,\u)\in\Q\times\U$,
the map $i^{\T}_{\q,\u}$ can be naturally split in two components.
Indeed,  for
   $(\p,\wp)\in T_{\q}^{*}\Q\times T_{\u}^{*}\U$, we can write
\bel{prodsplit}
   i^{\T}_{{\q},\u}(\p,\wp) ~=~ \left(\Big(i^{\T}_{{\q},\u}
   \Big)^{\Q}(\p,\wp),
   \Big(i^{\T}_{{\q},\u}\Big)^{\U}(\p,\wp)\right)~\in~   T_{\q}\Q
   \times T_{\u}\U\,.
   \eeq

\subsection{ Foliation structure and adapted coordinates}
The product structure of $\M=\Q\times\U$ induces  a foliation structure,
where the set of leaves is
\bel{foliation}{\Lambda} =\left\{ \Q\times \{\u\}\quad \u\in\U\right\}\,.
\eeq
       For every $(\q,\u) \in \Q\times \U$, we denote by
${\Lambda}(\q,\u) \doteq \Q\times \{\u\}$ the {\it leaf } through $(\q,\u)$.
   Let us consider the corresponding distribution\footnote{In
   our context, the term ``distribution"
    is meant in the sense of differential geometry,
namely, a fiber sub-bundle of the tangent bundle
   $T(\Q\times\U)$.} $\Delta$,
   whose fibers are given by
   $$
   \Delta_{(\q,\u)} = T_{\q}\Q \times\{0\}.$$
   In our analysis,
a very important role will also be played by the {\it orthogonal} distribution
\bel{orthbundle}
   \Delta^\perp_{(\q,\u)} = \Big\{ Y\in T_{\q}\Q\times T_{\u}\U\quad|\quad
   \g(\q,\u)(Y,X) = 0 \quad \forall X\in \Delta_{(\q,\u)}\Big\},
   \eeq
also called the {\it orthogonal bundle}, for short.

 In connection with the foliation $\Lambda$ at (\ref{foliation}),
 we say that a system of coordinates $(\tilde q,\tilde u)$
 is
 {\em $\Lambda$-adapted} if the sets $\{ \tilde u = \,\hbox{constant}\}$
 locally coincide with the leaves
 of the foliation.
 Of course, the local product coordinates
 $(q,u)$ are $\Lambda$-adapted.
More generally, if $(\tilde q,\tilde u)$
are  $\Lambda$-adapted, then every system of coordinates
$( \hat q , \hat u)$ obtained  from $(\tilde q,\tilde u)$
by means of a local diffeomorphism of the form
 \bel{adaptcoord}
 \hat q =  \hat q(\tilde q,\tilde u)\qquad  \hat u=  \hat u
 (\tilde u).
\eeq
is  $\Lambda$-adapted as well.

\vskip 1truecm \subsection{ Admissible input-output pairs }

Consider a control function $t\mapsto \u(t)\in \U$. In this section
we  characterize the corresponding output
 $t\mapsto (\q(t),\,  \p(t))$ as the solution of a certain Cauchy problem.
In the following  section, we then show that our definition
is consistent with the mechanical model, where the control is implemented in
terms of frictionless constraints.

 For every $(\u,\w)\in T\U\,$, let us define the map
$\T^{\u,\w}:T\Q\mapsto\rr$ by setting
\beq\label{kinetic-c} \T^{\u,\w}(\q,\v) \doteq \T(\q,\u,\v,\w)\,, \eeq
 for all $(\q,\v)\in T\Q$.
This map can be regarded as the kinetic energy of the reduced system
when the control takes the value $\u$, with $\dot \u=\w$.

 Let $I\subset\rr$ be an interval, and let
  $\u:I\mapsto\U$ be an absolutely continuous  control.
The (time-dependent) kinetic energy of the reduced system
 on $\Q$,
corresponding to the control $\u(\cdot)$ is described, almost every $t\in I$
and for all $(\q,\v)\in T\Q$, by
$$(t,\q,\v)\mapsto \T^{\u(t),\dot\u(t)}\big({\bf q},{\v }\big) .
$$
The corresponding (time-dependent) Hamiltonian on $T^*\Q$ is
$$
(t,\q,\p)\mapsto\H^{\u(t),\dot\u(t)}({{\q}},{\P})\,,
$$
where \bel{hamred} \H^{\u,\w}(\q,\p) \doteq \sup_{{\v }\in T_{{\q} }{\Q}
} \Big\{\langle{\p}, {\v }\rangle - \T^{\u,\w}({\q},\v)
\Big\} .\eeq
Since $\g$ is positive definite, for every
$(\q,\u)\in \Q\times\U$ and every $\p\in T_{\q}^*\Q$, the
 affine function
$$
\y\mapsto \Big(i^{\T}_{{\q},\u}\Big)^{\U}(\p,\y)~\in~ T_{\q}\Q
$$ in (\ref{prodsplit}) is
invertible. Its inverse will be denoted by
$$\w \mapsto
{\wp}_{(\q,\u,\p)}(\w).$$

 Let  $(q,u)$ be  $\Lambda$-adapted coordinates, and let
 $(q,u,p,{\pi})$ are the corresponding bundle coordinates.
 Let $(F_i,F_{N+\alpha})$ be the components of the force $\bf F$, so that
\bel{forzacomp0}
{\bf F} = F_i\frac{\partial}{\partial p_i} +
F_{N+\alpha}\frac{\partial}{\partial {\pi}_{\alpha}}\,. \eeq
Recalling the dimensions of the manifolds $\Q$ and $\U$, we here have
$i=1,\ldots, N$ and $\alpha = 1,\ldots,M$.
The Einstein summation convention is always used.
In addition, we set
\bel{forzaQ}
F_{\Q} \doteq F_i\frac{\partial}{\partial p_i}\,.
\eeq
Notice that $F_{\Q} $ is independent of the chosen
$\Lambda$-adapted system of coordinates.
For every
$(\u,\w)\in T\U$ and $i=1,\dots,N$, we also define
\bel{forza1}
F^{\u,\w}_i(t,\q,\p) \doteq
 F_i\big(t,\q,\u,\p,{\wp}_{(\q,\u,\p)}(\w)\big)
 \eeq
 and
 \bel{forza}
 {\bf F}_\Q^{\u,\w}(t,\q,\p) \doteq
 F_\Q\big(t,\q,\u,\p,{\wp}_{(\q,\u,\p)}(\w)\big) =
 F_i\big(t,\q,\u,\p,{\wp}_{(\q,\u,\p)}(\w)\big)
 \frac{\partial}{\partial p_i} \,.
\eeq
\vskip 3truemm
\begin{remark}
 {\rm Despite (\ref{kinetic-c}),
 in general one has
$$
\H^{\u,\w}\left(\q,\p\right) \neq \H\left(\q,\u,\p,
{\wp}_{(\q,\u,\p)}(\w)\right)\,.
$$
The actual relation between these two functions
will be illustrated in Lemma \ref{convex}.
}\end{remark}

 \begin{definition}\label{admissible}
     Let $I\subset \rr$ be a time interval. Let
     $$
     \u:I\mapsto\U \qquad \qquad (\q,\,
     \p):I\mapsto\T^*\Q
     $$ be  absolutely continuous maps.
      We say that  $\Big(\u(\cdot)\,,~
     (\q,\,
     \p)(\cdot)\Big)$
     is an {\em  admissible
     input-output pair}
    if $(\q,\,\p)$
is a Carath\'eodory solution  of the control equation of motion
\beq\label{eqham}
\frac{d}{dt}\Big(\q(t)\,,~
    \p(t)\Big) =
X_{\H^{\u(t),\dot\u(t)}}\Big(\q(t),\p(t)\Big) +
{\bf F}_\Q^{\u(t),\dot\u(t)}\Big(\q(t),\p(t)\Big)\, .\eeq
Here $X_{\H^{\u,\dot \u}}$ denotes the Hamiltonian
 vector field corresponding to
$\H^{\u,\dot \u}$, with respect to the symplectic structure
on  $T^*\Q$.
\end{definition}

We recall that a Carath\'eodory solutions of an ODE~ $\dot x= f(t,x)$ is an
absolutely continuous function $t\mapsto x(t)$ that satisfies the differential
equation at a.e.~time $t$.
Given an initial data
\beq\label{ic}\left(\q(\bar t), \p(\bar t)\right) =
\left(\bar{\q} , \bar{\bf p}\right)\,,
\eeq
and an absolutely continuous  control function $t\mapsto \u(t)$,
the existence and uniqueness of a corresponding admissible output
$(\q(\cdot), \,\p(\cdot))$ can be obtained from standard  ODE theory.


Depending on the geometrical properties
of the metric $\g$, the regularity assumptions
 on the input $\u$ and the output $(\q,\p)$
can be considerably weakened.
 This fact, discussed later on in the paper, is
 essential for both optimization and stabilization
 purposes.

    \subsection{ Realization of controls as frictionless constraints.} To motivate the previous notion
of input-output pair, we need to recall the notion
 of {\it frictionless constraint reaction} in the Hamiltonian framework.

    Let $Pr_1:T^*(\Q)\times T^*(\U)\to
     T^*(\Q)$ denote the canonical projection on the first factor, and, for every
      $((\q,\p),(\u,\y))\in T^*(\Q)\times T^*(\U)$, let
     $D(Pr_1)({\q},\p,\u,\y)$ denote its derivative.

 Let us consider the subspace of vertical vectors
 $$
 R_{(\q,\p),(\u,\y)}^{\Q} \doteq  \left(V_{(\q,\p)}(T^*\Q) \times
  V_{(\u,\y)}(T^*\U)\right) \bigcap \ker\big(D(Pr_1)({\q},\p,\u,\y)\big).
  $$
  It is straightforward to verify that
  $$
   R_{(\q,\p),(\u,\y)}^{\Q} = \{0\}\times
  V_{(\u,\y)}(T^*\U).
  $$
   \begin{definition}\label{friction}
  The subspace $R_{(\q,\p),(\u,\y)}^{\Q}$
  will be called the {\em subspace of $\Q$-frictionless reactions
  at $((\q,\p),(\u,\y))$.}
The corresponding  vector  bundle based on $T^*{\Q}\times T^*{\U}$
  will be called
  {\em the vector bundle of $\Q$-frictionless reactions}.
   \end{definition}
  \begin{remark} {\rm In terms of the canonical form $\omega_{\M}$ on $\M=\Q\times\U$, $R_{(\q,\u),(\p,\y)}^{\Q}$ can be characterized as the subspace
  of $V_{(\q,\p)}(T^*\Q) \times
  V_{(\u,\y)}(T^*\U)$ made
   of
  those vectors $\Phi$ such that
  \bel{liscio}
  \omega_\M(\Phi,V) = 0 \qquad \forall {\bf V}\in T_{(\q,\p)}(T^*Q)\times \{0\}
  \eeq
We remind that, in  the language of symplectic geometry,  $\Phi$ and $\bf V$ are also said  {\em anti-orthogonal}.}
\end{remark}
\begin{remark}{\rm We are here  regarding the constraint reactions as a set-valued force, described by the multifunction
  $$
{(\q,\p),(\u,\y)}\mapsto R_{(\q,\p),(\u,\y)} .
$$  To check that this definition coincides with the usual one
it is sufficient to notice that   if $\Phi \in R_{\q,\p,\u,\y}^{\Q}$ and  $\Phi=\sum_{r=1}^{N+M}{\Phi}_r
   \frac{\partial}{\partial P_\RR}$ is its local expression, then (\ref{liscio}) is equivalent to
   $$ \left\langle \left(\Phi_1,\dots,\Phi_{N},\Phi_{N+1},\dots, \Phi_{N+M}\right)\, , \, \left(v_1,\dots, v_N,0,\dots,0\right)
   \right\rangle = 0\qquad \hbox{for all} \quad v_1,\dots,v_N\in \rr.
   $$
  Of course, this holds if and only if $
   {\Phi}_i = 0$ for all $  i=1,\dots,N
   $.}
\end{remark}

Definition \ref{admissible} is justified by Theorem \ref{equivalence} below.  Let $I\subset\R$ be an interval, and, for every $(t,\q,\u,\p,\y)\in I\times T^*Q\times T^*\U\to T^*\Q$, let us set
   $$
   X^{\Q}_{\H}(t,\q,\u,\p,\y)\doteq
D(Pr_1)\cdot X_{\H}(t,\q,\u,\p,\y).$$
Notice that,   according to (\ref{forza}),
$$
 F_{\Q}(t,\q,\u,\p,\y)\doteq
D(Pr_1)\cdot F(t,\q,\u,\p,\y).
$$
 \vskip 1truecm
\begin{thm}\label{equivalence}  Consider\,\,\,\, a time\,\,\, interval \,\,\, $I\subset\rr$
and \,\,\,let the maps
     $\u:I\mapsto\U\,$, ~$(\q,\,\p):I\mapsto T^*\Q$ be twice continuously
differentiable. Then the following conditions are equivalent:
       \vskip5truemm
      {\bf(i)}~ $\Big(\u(\cdot)\,,~(\q(\cdot),
     \p(\cdot))\Big)$
     is an {\it  admissible
     input-output pair}, that is, $(\q,\,\p)$
verifies
\beq\label{eqham1}
\frac{d}{dt}\Big(\q(t)\,,~
    \p(t)\Big) =
X_{\H^{\u(t),\dot\u(t)}}\Big(\q(t),\p(t)\Big) +
{\bf F}_\Q^{\u(t),\dot\u(t)}\Big(\q(t),\p(t)\Big)\, .\eeq
     \vskip5truemm
      {\bf(ii)}  The path
  $\Big(\q(\cdot) ,
     \p(\cdot)\Big)$
is an integral curve of the control system
\begin{equation}\label{44.16}
\frac{d}{dt}\Big(\q(t) ,
   \p(t)\Big) =
    \Big[X^{\Q}_{\H}(t,\q(t),\u(t),\p(t),\y) +
     F_{\Q}(t,\q(t),\u(t),\p(t),\y)\Big]_{|_{
    \y = {\wp}_{(\q(t),\u(t),\p(t))}(\dot\u(t))}}\,.
\end{equation}
\vskip5truemm
     {\bf(iii)} ~ There exist selections
    $$
    t\mapsto \y(t)\in T_{u(t)}^*(\U)
    \qquad\qquad t\mapsto r(t)\in R^\Q_{\q(t),\p(t),\u(t),\y(t)}
    $$
    such that, for all $t\in I$, one has
     \beq\label{motionincl}
 \frac{d}{dt}\Big(\q(t),\, \u(t),\,
    \p(t)\,,\y(t)\Big)  = X_{\mathcal H}
 + {\bf F}
  + r(t)\,.
 \eeq
 The map $r(\cdot)$ in {\rm (\ref{motionincl}) } is called the {\em constraint reaction
 corresponding to the motion $(\q,\p,\u,\y)(\cdot)$}.
\end{thm}

{\bf Proof. }
Assume that condition {\bf(iii)} holds.
      In particular, for every $t\in I$,
    one has
    $$
    \dot\u(t) =
     \Big(i^{\T}_{{\q(t)},\u(t)}\Big)^{\U}(\p(t),\y(t)),
    $$
    which implies
    $$
    \y(t) = {\wp}_{(\q(t),\u(t),\p(t))}(\dot\u(t)).
    $$
    Therefore
    $$\begin{array}{rl}
    \frac{d}{dt}\Big(\q(t) ,
    \p(t)\Big) & =~\displaystyle
     D(Pr_1)\cdot \frac{d}{dt}\Big(\q(t) ,
    \u(t),
    \p(t), \y(t)\Big)~ =~
     D(Pr_1 )\cdot \Big( X_{\mathcal H} + F +
     r(t)\Big) \\\, \\
    &= ~X_{\mathcal H}^{\Q}(\q(t),\u(t),\p(t),\y(t))
     + F_{\Q}(\q(t),\u(t),\p(t),\y(t)) \\ \, \\
&=~\Big[X^{\Q}_{\H}(\q(t),\u(t),\p(t),\y) +
     F_{\Q}(\q(t),\u(t),\p(t),\y)\Big]_{|_{
    \y ={\wp}_{(\q(t),\u(t),\p(t))}(\dot\u(t))}}\,.
    \end{array}
    $$
    Hence {\bf (ii)} holds as well.

Conversely, let us show that {\bf (ii)} implies condition {\bf (iii)} .
If $\Big(\q(\cdot) ,
     \p(\cdot)\Big)$  is a solution of (\ref{44.16}),
 then, setting
$$
\y(t)= {\wp}_{(\q(t),\u(t),\p(t))}(\dot\u(t))\,,$$ one has
$$
\dot\u(t) =
     \Big(i^{\T}_{{\q(t)},\u(t)}\Big)^{\U}(\p(t),\y(t))
    $$
    which coincides with the equation for the variable $\u$ in
    (\ref{motionincl}).
Therefore, condition {\bf (iii)}  is
satisfied provided that we define
$$
r(t) \doteq \dot\y(t) - X_{\H}\big(\q(t),u(t),\p(t),\y(t)\big) -
F\big(\q(t),u(t),\p(t),\y(t)\big).
$$

     \vskip8truemm
     In order to prove that conditions {\bf (i)}
      and  {\bf (ii)} are equivalent,
     it is sufficient to prove that the right-hand sides
      of the involved equations do coincide.  Actually,
      recalling the definition
     of $F_{\Q}^{\u,\w}$ at  (\ref{forza}), for all $\big((\q,\p),(\u,\w)\big)\in
     T^*\Q\times T^*\U$, one has
     $$
     \Big[F_{\Q}(t,\q,\u,\p,\y)\Big]_{|_{
    \y = {\wp}_{(\q,\u,\p)}(\w)}}  =
F^{\u,\w}\Big(\q,\p\Big).
$$
Hence the proof is concluded as soon  as one proves  the identity
\beq
     \Big[X^{\Q}_{\H}(\q,\u,\p,\y)\Big]_{|_{
    \y = {\wp}_{(\q,\u,\p)}(\w)}}  = X_{\H^{\u,\w}}\Big(\q,\p\Big)
\eeq This will be done  in Section \ref{loc}  ---see  Lemma \ref{convex} below--- by use of  of local coordinates.

\endproof


\section{The control equation in local coordinates}\label{loc}
\setcounter{equation}{0}

Consider a  $\Lambda$-adapted coordinate chart
$(q,u)$ defined on an open set $U$, and let
$\Big((q,u),(p,w)\Big)$ be the
corresponding coordinates on the fiber bundle
 $$\bigcup_{(\q,\u)\in U}\{(\q,\u)\}\times (T_\q^*\Q\times T_\u\U).$$
Let
$G=(g_{r,s})_{r,s=1,\dots,N+M}$ be the matrix representing the
kinetic metric $\g$, and let $G^{-1}=(g^{r,r})_{r,s=1,\dots,N+M}$ denote
its inverse.  In the following, we consider the sub-matrices
$$G_1\doteq(g_{i,j})\,,\qquad
  G_2\doteq \left(g_{N+\alpha,N+\beta}\right)\,,\qquad
 (G^{-1})_2\doteq
 \left(g^{N+\alpha,N+\beta}\right)\,,
 $$
$$ G_{12}\doteq
 \left(g_{i,N+\alpha}\right)\,,\qquad
 (G^{-1})_{12}\doteq \left(g^{i,N+\alpha}\right)\,,
$$
with the convention that the Latin
 indices $i,j$ run from $1$ to $N$, while the Greek indices
 $\alpha,\beta$ run from $1$ to $M$.
For convenience, we also define
\bel{matrixdef}
A=(a^{i,j}) \doteq (G_1)^{-1},\qquad E=(e_{\alpha,\beta})
\doteq ((G^{-1})_2)^{-1},
\qquad K=(k^i_{N+\alpha}) \doteq  (G^{-1})_{12}E\,.
\eeq

Let $\u(\cdot):I\mapsto\U$ be twice continuously differentiable, and let
$(\q ,\p):I\mapsto T^*\Q$
be a curve such that $\Big(\u(\cdot)\,,\,(\q,\p)(\cdot)\Big)$
is an admissible input-output
pair for the control equation of motion. By possibly restricting
the size of the interval $I$, we can assume that
$\Big(\q(t),\u(t)\Big)$ remains inside the domain of the single chart
$(q,u)$ for every $t\in I$.

\begin{thm}\label{eqonNloc} Given an admissible
input-output pair, the corresponding coordinate maps
$t\mapsto\left(u(t)\,,\,\pmatrix{ q(t)
\cr p^\dagger(t)}\right)$ satisfy  the  differential equation
\beq\label{semihamex} \begin{array}{c}
\left(\begin{array}{c}\dot q\\\,\\\dot p\end{array}\right) =
 \left(\begin{array}{c} Ap \\\,\\
-\frac{1}{2}p^\dagger\frac{\partial A}{\partial
q}p\end{array}\right)
+\left(\begin{array}{c} { K\dot u}\\\,\\
-p^\dagger\frac{\partial { K}}{\partial q}\dot{u} \end{array}\right)  +
  \left(\begin{array}{c} {0}\\\,\\
\frac{1}{2}\dot u^{{\dag}}\frac{\partial E}{\partial
q}\dot{u}\end{array}\right)  +  \left(\begin{array}{c}0\\\,\\
F_\Q^{{u}(\cdot),\dot{u}(\cdot)}
\end{array}\right) \,.\end{array}
\eeq
In addition, recalling (\ref{forza1})-(\ref{forza}),
\bel{13}
F_\Q^{{u}(\cdot),\dot{u}(\cdot)}\doteq \Big(F_1^{{u}(\cdot),\dot{u}(\cdot)},
\dots,F_N^{{u}(\cdot),\dot{u}(\cdot)}\Big) .\eeq
\end{thm}

For convenience, in (\ref{semihamex}) we write all vectors as column vectors,
while the superscript $^{{\dag}}$ denotes transposition.  Componentwise, this means:
\beq\label{semihamex1}
\left\{ \begin{array}{lll} \dot q^i = &\,\,\,\,a^{i,j}p_j\,\,\,\,\,\,
\,\,\,\,\,\,\,\,\, +&k^i_{N+\alpha} \dot u^{\alpha}\,,\\\,\\
\dot p_i = &-\frac{1}{2}\frac{\partial a^{\ell,j}}{\partial q^i} p_\ell p_j\,\,\,\, -\,\,
&\frac{\partial k_\alpha^j}{\partial q^i}p_j \dot u^\alpha
+ \frac{1}{2}\frac{\partial e_{\alpha,\beta}}{\partial q^i} \dot u^\alpha
\dot u^\beta + F_i^{{u}(\cdot),\dot{u}(\cdot)}.
\end{array}\right.
\eeq
(where $\ell$ runs from $1$ to $N$).

{\bf Proof. }
Let $\H^{u,w}$ be the local coordinate representation of the Hamiltonian
$\H^{\u,\w}$.   Then the map $(q ,p)(\cdot)$ satisfies the system
\beq\label{semiham} \left\{\begin{array}{l}\dot q^i(t) =
\frac{\partial \H^{u(t),\dot u(t)}}{\partial p_i}(q(t),p(t))\,,
\\\,\\
\dot p_i(t) =-\frac{\partial \H^{u(t),\dot u(t)}}{\partial
q^i}(q(t),p(t)) + F_i^{{u}(t),\dot{u}(t)} (t,q(t),p(t))\,.
\end{array}\right.
\eeq

Using the  expression of the kinetic energy
$$
\T^{u,w} = \frac{1}{2} \dot q^{{\dag}} G_1 \dot q +
\dot q^{{\dag}} G_{12}  w + \frac{1}{2} w^{{\dag}} G_2 w
$$
 in terms of the given coordinates, we obtain
 \bel{Hcoord}\H^{u,w}(q,p) =
\H^{u,w}_2(q,p) + \H^{u,w}_1(q,p) + \H^{u,w}_0(q,p)\,,\eeq where
$$
\H^{u,w}_2(q,p) = \frac{1}{2}p^\dagger{ {(G_1)^{-1}}}p\,,\qquad\quad
\H^{u,w}_1(q,p) = -p^\dagger{ {(G_1)^{-1}}G_{12}}w
\,,
$$
$$
\H^{u,w}_0(q,p) = \frac{1}{2} w^{{\dag}} {
G_{12}^{{\dag}}{(G_1)^{-1}}G_{12}}w - \frac{1}{2}w^{{\dag}}{ G_2}w\,. $$
The theorem can thus be proved by checking that
\beq\label{identici}
\H^{u,w}(q,u) = \frac{1}{2}p^\dagger{
A}p + p^\dagger{ K} w - \frac{1}{2}
 w^{{\dag}}{ E} w\,.
 \eeq
{}From the identity
$$
{ G_1(G^{-1})_{12}} + { G_{12}(G^{-1})_2} = 0
$$
one obtains
$$ 0= { {G_1^{-1}}G_1(G^{-1})_{12}E} + {
{G_1^{-1}}G_{12}(G^{-1})_2E},
$$
 which implies
\beq\label{primamat}
 \Big(K =\Big)\,\,\,\,\,   { (G^{-1})_{12}E} = {-{(G_1)^{-1}}G_{12} }
 =  {-{A}G_{12} }  \,. \eeq
 Moreover,  from  (\ref{primamat}) and the identity
 $$
 { \bf I}_M =   G_{12}^{{\dag}}{ (G^{-1})_{12}} + { G_2(G^{-1})_2}\,,
 $$
one gets
\beq\label{secondamat}
{ E} = G_{12}^{{\dag}}{ (G^{-1})_{12}E} + { G_2(G^{-1})_2E} = - G_{12}^{{\dag}}{
AG_{12}} + { G_2}\,.
\eeq
Together, (\ref{primamat}) and (\ref{secondamat}) yield (\ref{identici}),
concluding the proof.
\endproof

  We conclude this section by proving  the following lemma, which was used in the proof of
  Theorem \ref{equivalence}.

 \begin{lma}\label{convex}  For
    all $\big((\q,\p),(\u,\w)\big)\in
     T^*\Q\times T\U$
     one has
     \beq\label{convexeq}
    \left[X^{\Q}_{\H}(\q,\u,\p,\y)\right]_{|_{
    \y = {\wp}_{(\q,\u,\p)}(\w)}} =   X_{\H^{\u,\w}}\Big(\q,\p\Big)
     \eeq
     \end{lma}

{\bf Proof.}  Since the thesis is local one can use coordinates to prove it.
Consider a  $\Lambda$-adapted coordinate chart
$(q,u)$  and let
$\Big((q,u),(p,w)\Big)$  $\Big((q,u),(p,\wp )\Big)$ be the
corresponding coordinates on the fiber bundles
 $ T^*\Q\times T\U$  and
  $T^*\Q\times T^*\U$, respectively.
  Let $(G^{-1})_{12}, A, E, K$
  have the same meaning as before.
  Moreover, let us set $((G^{-1})_1 \doteq (g^{i,j})_{i,j=1,\dots,N}$.
  We have to prove that
  \bel{f1}
  \frac{\partial \H^{u,w}}{\partial p} = \frac{\partial \H}{\partial p}_{|_{
    \wp = \wp_{(q,u,p)}(w)}}
    \eeq
    and
    \bel{f2}
 \frac{\partial \H^{u,w}}{\partial q} = \frac{\partial \H}{\partial q}_{|_{
    \wp = \wp_{(q,u,p)}(w)}},
    \eeq
    where
 $$
  \wp_{(q,u,p)}(w) \doteq w^{{\dag}}E + p^\dagger K.
  $$

In fact, by (\ref{Hcoord}) one has
 \bel{f3} \frac{\partial \H^{u,w}}{\partial p}  = p^\dagger A+Kw.
 \eeq
On the other hand, one can easily check that
\bel{f4}
 \frac{\partial \H}{\partial p}_{|_{
    \wp = \wp_{(q,u,p)}(w)}} = p^\dagger (G^{-1})_1 +  (G^{-1})_{12} \, (w^{{\dag}}E
    + p^\dagger K )= p^\dagger A +Kw,
     \eeq
    so (\ref{f1}) is proved.

To prove (\ref{f2}), let us consider the well-known identities
      $$
    \frac{\partial{\H}}{\partial q}
     (q,u,p,{\wp}) = - \frac{\partial{\T}}{\partial q}
     (q,u,v,w)_{|_{(v,w) = \left(\Big(i^{{\H}}_{{\q},\u}
 \Big)^{\Q}(p,{\wp}),\Big(i^{{\H}}_{{\q},\u}
 \Big)^{\U}(p,{\wp})\right)}}$$
 $$
 \frac{\partial{\H^{u,w}}}{\partial q}
     (q,p) = - \frac{\partial{\T^{u,w}}}{\partial q}
     (q, v)_{|_{v= \frac{\partial{\H^{u,w}}}{\partial p}(q,p)}}
     $$
 Observe that the map
$p\mapsto  \frac{\partial{\H^{u,w}}}{\partial p}(q,p)$
is the inverse of the map $v\mapsto  \frac{\partial {\T}^{u,w}}{\partial v}(q,v)$.
    Hence, letting
    $$
    \check v \doteq\Big(i^{{\H}}_{{\q},\u}
 \Big)^{\Q}\left(p,\Big(i^{\T}_{{\q},\u}
 \Big)^{\U}\left( \frac{\partial{\H^{u,w}}}{\partial p},w\right)\right)
 $$
 $$
 \check w \doteq  \Big(i^{{\H}}_{{\q},\u}
 \Big)^{\U}\left(p,\Big(i^{\T}_{{\q},\u}
 \Big)^{\U}\left( \frac{\partial{\H^{u,w}}}{\partial p},w\right)\right)\,,
 $$
 one obtains
 $$
      \frac{\partial{\H}}{\partial q}
     (q,u,p,{\wp})_{|_{{\wp}={\hat \wp}(q,u,p,w)}} =
     \frac{\partial{\H}}{\partial q}
 (q,u,p,{\wp})_{|_{{\wp}= \Big(i^{\T}_{{\q},\u}
 \Big)^{\U}\left( \frac{\partial{\H^{u,w}}}{\partial p},w\right)}} = $$
 $$
 - \frac{\partial{\T}}{\partial q}
     (q,u,v,w)_{|_{(v,w) = (\check v,\check w) }}=
  - \frac{\partial{\T}}{\partial q}
     (q,u,v,w)_{|_{v= \frac{\partial{\H^{u,w}}}{\partial p}(q,p)}} =
      \frac{\partial{\H^{u,w}}}{\partial q}
     (q,p)\,\,\,,
     $$
so (\ref{f2}) is proved as well.
\endproof

\begin{remark} {\rm As a matter of fact, one could show that  Lemma \ref{convex} is valid under the weaker assumption that $\T$ is strictly convex in the velocity $(\v,\w)$ --- see \cite{RLect}.} \end{remark}

\section{ The Riemannian structure and and the
input-output map}\label{ap1}\setcounter{equation}{0}

The presence of the derivative
$\dot\u$ in the dynamic equations (\ref{semihamex})
depends on the Riemannian metric $g$ defining the kinetic energy
and on the foliation
$\Lambda$ at (\ref{foliation}).
In this section we review the main results in this direction.
To simplify the discussion, throughout this section we shall assume that
the additional  forces ${\bf F}$ vanish identically, so that
in (\ref{13}) one has
$$F_\Q^{{u}(\cdot),\dot{u}(\cdot)}\equiv 0\,.$$

The following definitions were introduced in \cite{AB1}.

     \begin{definition} A local,
     $\Lambda$-adapted, system of coordinates $(q,u)$ on
$\Q\times\U$ is called
     {\em $N$-fit for hyperimpulses } if,
     for every  differentiable control function $\u(\cdot)$,
     the right-hand side

     of the corresponding  equation of motion (\ref{semihamex1})
     does not contains any quadratic term
     in the variable $\dot u$.

     A local,   $\Lambda$-adapted, system of coordinates $(q,u)$
     on $\Q\times\U$ is called
     {\em strongly  $N$-fit for hyperimpulses } if,
     for every  differentiable control function $\u(\cdot)$,
     right-hand side of (\ref{semihamex1})
     is independent of the variable $\dot u$.

      Moreover, we shall call {\em generic} any local ,
      $\Lambda$-adapted, system of coordinates
      $(q,u)$ which is not $N$-fit for hyperimpulses.
      \end{definition}

   \begin{remark}  {\rm The denomination ``{\em $N$-fit for hyperimpulses}"
   for a system of coordinates   $(q,u)$ refers to the fact
  that,
if the dependence on $\dot u$ is only linear,
one can then construct solutions $\Big({q}(\cdot),p(\cdot)\Big)$
also for discontinuous controls $u(\cdot)$.
In general, a jump in $u(\cdot)$ will produce
a discontinuity  in both components
${q}(\cdot)$ and $p(\cdot)$.  For this reason we call it
a  {\em hyperimpulse}, as opposite to  {\em impulse}, which
can cause a discontinuity in the component  $p(\cdot)$ only.
}\end{remark}

A first characterization of $N$-fit coordinates  was derived in
\cite{AB1}.
It is important to observe that the property of being $N$-fit depends only on
the metric $g$ and on the foliation $\Lambda$, not on particular the system of
$\Lambda$-adapted coordinates.  This leads to




      \begin{definition}\label{fitfol}{\rm \cite{Rampazzo2}} The foliation $ {\Lambda}$ is called
{\em $N$-fit for hyperimpulses } if there exists an atlas of
$\Lambda$-adapted charts that are also $N$-fit for hyperimpulses.
In this case, {\em all}
$\Lambda$-adapted charts are then  $N$-fit for hyperimpulses.

The foliation $ {\Lambda}$ is called {\em strongly $N$-fit for
hyper-impulses } if there exists an atlas of $\Lambda$-adapted
charts which are strongly $N$-fit for hyper-impulses.

Moreover, the foliation ${\Lambda}$ will be
       called {\em generic  } if it is not $N$-fit for hyper-impulses.
       \end{definition}

The paper \cite{Rampazzo2}
established the connection between the $N$-fitness of the foliation
$\Lambda$
and the bundle-like property of the metric, introduced in
\cite{Re1, Re2}.
We recall here the main definitions and results.

       \begin{definition} The metric $\g$ is {\em bundle-like}
       with respect to the foliation $\Lambda$
      if,  for every $\Lambda$- adapted chart,
      it has a local representation of the form
       $$
       \sum_{i,j=1}^N g_{i,j}(q,u)\omega^i\otimes\omega^j +
       \sum_{\alpha,\beta=1}^M g_{N+\alpha,N+\beta}(q,u) dc^{\alpha}\otimes
       dc^{\beta}\,,
       $$
       where $\omega^1,\dots,\omega^N$ are linearly independent $1$-forms
       such that, for each $(\q,\u)\in \Q\times \U$ in the domain of the chart,
       one has

       (i) $\Big(\omega^1(\q,\u),\dots,\omega^N(\q,\u),dc^1(\u),
       \dots,dc^M(\u)\Big)$ is a basis of  the cotangent space
       $T^*_{\q}\Q\times T^*_{\u}\U$;

       (ii) $\Big\langle \omega^i(\q,\u), Y\Big\rangle = 0$,
       for every $Y\in \Delta^\perp_{(\q,\u)}$.
       \end{definition}

We recall that $\Delta^\perp_{(\q,\u)}$ is  the orthogonal bundle,
defined at (\ref{orthbundle}).
If $g$ is bundle-like with respect to the foliation $\Lambda$,
the latter is also called a {\em Riemannian foliation},
because in this case
a Riemannian structure can be well defined also on the quotient space.
In order to state the next theorem, we recall  the notion
 of {\em completely integrable}
  distribution. (Here, by {\it distribution} we mean a map that maps

  \begin{definition} Let $\M$ be a  manifold of dimension
  $d$,
  and let $\Gamma$ be a distribution on
  $\M$ of dimension $ N\leq  d$.  (I.e., for  every  $y\in M$, $\Gamma(y)$ is a subspace of dimension $N$ of $T_y\M$.) We say that the distribution $\Gamma$
  is {\em completely  integrable} if,
  for every $\m\in\M$, there exists  a neighborhood $U$ of $\m$ and a local
system of coordinates $(x,z) = (x^1,\dots,x^N,z^{N+1},\dots,z^d)$
such that at each point $\m\in U$ one has
  $$
  \Gamma_{\m} = span\left\{ \frac{\partial}{\partial x_i}\qquad
i=1,\dots,N\right\}\,.
  $$
  \end{definition}

\begin{thm}\label{mfit-th}  On the product space $\M=\Q\times\U$,
consider the natural foliation $\Lambda$ as in (\ref{foliation}).
The following statements are equivalent:

i) The foliation $ {\Lambda}$ is $N$-fit for hyper-impulses.

ii) The metric $\g$ is bundle-like with respect to the foliation
 ${\Lambda}$, i.e., the foliation $ {\Lambda}$ is Riemannian.

 iii) For any $\u,\bar\u\in \U$ the map
 $d_{\u,\bar\u}(\cdot): \Q\mapsto \rr$ defined by
 $$
 d_{\u}(\q) \doteq dist\,\,\Big((\q,\u), \Q\times\{\bar\u\}\Big)
 $$ is constant. In other words,
 leaves remain at the same distance from each other.

 iv) If $t\mapsto(\q(t),\,\u(t))$  is any geodesic curve with
respect to the metric $\g$, and if $(\dot \q(\tau),\dot \u(\tau))\in
\Delta^\perp_{(\q(\tau),\u(\tau))}$ at some time $\tau$, then
$(\dot \q(t),\dot \u(t))\in
\Delta^\perp_{(\q(t),\u(t))}$ for all $t$.   In other words,
if a geodesic crosses perpendicularly one of the leaves, then it
crosses perpendicularly also every other leaf which it meets.




v) If $(q,u)$ is a $\Lambda$-adapted system of coordinates, then
\bel{N-fitcondition}
\frac{\partial g^{N+\alpha ,N+\beta}}{\partial q_i} = 0 \qquad \qquad
i=1,\dots,N,\quad
\alpha,\beta=1,\dots,M\,,
\eeq
where $G^{-1} =(g^{r,s})$ denotes the inverse of the matrix $G=(g_{r,s})$
representing the metric $\g$ in the coordinates $(q,u)$.

 \end{thm}

 Indeed, the equivalence of i) and ii) is a trivial consequence of the definitions
 of bundle-like metric and of  $N$-fit system of coordinates.
 The equivalence of ii), iii), and iv), is a classical result on
 bundle-like metrics \cite{Re1}.
 Moreover, by (\ref{semihamex1}), the foliation is fit for jumps
 if and only if $\partial e_{\alpha,\beta}/\partial q^i\equiv 0$.
 Recalling that the matrix $(e_{\alpha,\beta})$ is the inverse of
 $(G^{-1})_2=(g^{N+\alpha, N+\beta})$, we conclude  that i) is equivalent to v).

\begin{thm}\label{st-mfit-th} The following statements are equivalent:

i) The foliation $ {\Lambda}$ is strongly $N$-fit for hyperimpulses .

ii) The foliation $ {\Lambda}$ is $N$-fit for hyperimpulses and the
orthogonal bundle $\Delta^\perp_{(\q,\u)}$ in (\ref{orthbundle})
is integrable.

iii) There is an atlas such that, for every  chart $(q,u)$, one has
$$
\frac{\partial g^{N+\alpha, N+\beta}}{\partial q_i} = 0,
\qquad  g^{i,N+\alpha} = 0\qquad \forall i=1,\dots,N,\quad
\alpha,\beta=1,\dots,M.
$$



  \end{thm}
Indeed the equivalence of i) and ii),
formulated in terms of  Riemannian foliations,
was proved in  \cite{Re1}.
The equivalence between i) and iii)
follows from  (\ref{semihamex1}).
See again \cite{Re1} and \cite{Rampazzo2} for details.

\section{ The quadratic term in the control equation and the
orthogonal curvature
of  the foliation $\Lambda$}\label{ap2}\setcounter{equation}{0}

As we have seen in the previous section, the $N$-fitness for hyperimpulses
of a coordinate system $(q,u)$ can be characterized in terms of geodesics.
Indeed, the quadratic terms in the control equation
of motion (\ref{semihamex1})
are identically zero if and only if any geodesic which crosses
perpendicularly one leaf of the foliation $\Lambda$ also has
perpendicular intersection with every other
leaf it meets.

In the general case, however, the quadratic terms in
(\ref{semihamex1}) do not vanish.  We wish to give here a geometric
interpretation of these terms. This will again  be achieved
by looking at geodesics whose  tangent vector initially
lies in the orthogonal distribution
$\Delta^\perp$.

\subsection{$\U$-orthonormal coordinates}
 We shall make an essential use of a Proposition\ \ref{orthoC} below,
 which establishes the existence of a special kind  of $\Lambda$-adapted charts.
 To state it, let us set $(x^1,\dots, x^{N+M}) \doteq (q^1,\dots,q^N,u^1,
 \dots,u^M)$,  and, for every $h,k,r,s = 1,\dots,N+M$, let us consider the  functions
 $$
 \Gamma_{h,r,s} \doteq \frac{1}{2}\left(\frac{\partial g_{h,r}}{\partial x^s} + \frac{\partial g_{h,s}}{\partial x^r} -
 \frac{\partial g_{r,s}}{\partial x^h}\right)\qquad  \Gamma^k_{r,s} \doteq g^{kh}\left(\frac{\partial g^{h,r}}{\partial x^s} + \frac{\partial g^{h,s}}{\partial x^r} -
 \frac{\partial g^{r,s}}{\partial x^h} \right)
 $$
The $ \Gamma^k_{r,s} $ are the well-known Christoffel symbols.
 \begin{prop} \label{orthoC} Consider a point $(\bar\q,\bar\u)\in  \Q\times\U$
 and an orthonormal basis $\{J_1,\dots,J_M\}$ of $ \Delta^{\perp}(\bar\q,\bar \u)$.
 Then there exist
 $\Lambda$-adapted
 coordinates  $ (q,u)$, defined on a neighborhood of  $(\bar\q,\bar\u)$,  such that   calling
 $G=(g_{r,s})_{r,s=1,\ldots, N+M}\,$
 the corresponding kinetic
 matrix, one has

  {\rm (i)} the point $(\bar\q,\bar\u) $ has coordinates $(0,0)$;

{\rm (ii)}~ $g_{r,s} (0,0)= g^{r,s} (0,0)= \delta_{r,s}$
(the Kronecker symbol) for all
$r,s=1,\dots,N+M$;

{\rm (iii)} for every $w = (w_1,\dots,w_M)\in \R^M$, the geodesic
$ (q,u)_{w}(\cdot)$
issuing from $
(\q,\u)$ with velocity equal to $ w_1J_1+\dots w_MJ_M$ has  local representation
$(q,u)_{w}(t) = (0,\dots,0,tw_1,\dots,tw_M)\,.
$

Moreover, for all indices $i=1,\ldots,N$ and $\alpha,\beta,\gamma=1,\ldots,M$,
we have
\bel{Chr}
\Gamma_{i, N+\alpha,N+\beta}(0,0)  = \Gamma^i_{N+\alpha,N+\beta}(0,0) = 0 , \quad \Gamma_{N+\gamma, N+\alpha,N+\beta}(0,0) = \Gamma^{N+\gamma}_{N+\alpha,N+\beta}(0,0) =0
\eeq
In turn, this implies
\bel{equal1}
\frac{\partial g^{N+\beta,N+\gamma}}{\partial q^i}
\,(0,0) ~=~\frac{\partial g^{i, N+\beta}}{\partial u^\gamma}\,(0,0)
+\frac{\partial g^{i, N+\gamma}}{\partial u^\beta}\,(0,0)\,,
\eeq
\bel{equal2}
\frac{\partial g_{N+\alpha,N+\beta}}{\partial u^{\gamma}}(0,0)
~=~\frac{\partial g^{N+\alpha,N+\beta}}{\partial u^{\gamma}}(0,0)
~ = ~0\,.
\eeq

A chart with the above properties will be called {\em $\U$-orthonormal at
 $(\bar\q,\bar\u)$}.
\end{prop}

{\bf Proof.}  We start by considering  $\Lambda$-adapted coordinates
$(\hat q,\hat u)$, defined on a neighborhood of the point $(\bar\q,\bar \u)$,
 such that at the point $(\bar\q,\bar \u)$ one has
$(\hat q,\hat u) = (0,0)$ and
$$\frac{\partial}{\partial \hat u^\alpha} = J_\alpha\,,\qquad\qquad \alpha
=1,\ldots,M\,,$$
while
$$\left\{\frac{\partial}{\partial \hat q^i}\,;~~i=1,\ldots,N\right\}$$
is an orthonormal basis of the tangent space $T\Q$ at $(\bar\q,\bar \u)$,
w.r.t.~the metric~${\bf g}$.

To achieve the further property (iii), we need to modify these coordinates,
using the exponential map.
In the following, given a tangent vector
${\bf V}\in T_{(\bar \q,\bar\u)}(\Q\times\U) $, we denote
by
$\tau\mapsto \gamma_{\bf V}(\tau)$ the geodesic curve
starting from $(\bar\q,\bar\u)$ with velocity $\bf V$. In other words,
$$\gamma_{\bf V}(0)=(\bar\q,\bar\u)\,,
\qquad\qquad \frac{d\gamma_{\bf V}}{d\tau}(0)=
{\bf V}\,.$$
The exponential map is then defined by setting
$$
Exp_{(\bar\q,\bar\u)}({\bf V}) \doteq \gamma_{\bf V}(1).
$$
This is well defined for all vectors ${\bf V}$ in a neighborhood of the origin.

Denote by $\Big(\hat q,\hat u\Big)(\q,\u)$
the coordinates of a point  $(\q,\u)$ via the chart $(\hat q,\hat u)$.
We now define $(q,u)$ to be the new  coordinates of a point $(\q,\u)$
provided that
\bel{orthcoordef}
\Big(\hat q,\hat u\Big)(\q,\u)= \Big(\hat q,\hat u\Big)
\left(Exp_{(\bar\q,\bar\u)}\Big(\sum_{\alpha=1}^M
u^{\alpha}J_{\alpha}\Big)\right) + (q^1,\dots,q^N,0,\dots, 0).
\eeq
Notice that this is well defined, for all
$(\q,\u)$ in a neighborhood of $(\bar \q,\bar \u)$.
Indeed, the map $\rho:\rr^{N+M}\mapsto\rr^{N+M}$ defined by
\bel{rho}
\rho(q^1,\ldots,q^N, u^1,\ldots, u^M)= \Big(\hat q,\hat u\Big)
\left(Exp_{(\bar\q,\bar\u)}\Big(\sum_{\alpha=1}^M
u^{\alpha}J_{\alpha}\Big)\right) + (q^1,\dots,q^N,0,\dots, 0).
\eeq
maps the origin into itself. Moreover, by the properties
of the chart $(\hat q,\hat u)$,
the Jacobian matrix
$\partial \rho/\partial (q^1,\dots,q^N,u^1,\dots, u^M)$
at the origin coincides with  the identity matrix.
This already establishes the properties (i)-(ii).

By construction, $(q,u)(\q,\u)= (0,\ldots, 0, u^1,\ldots,u^M)$
if and only if
$$(\q,\u)=Exp_{(\bar\q,\bar\u)}\Big(\sum_{\alpha=1}^M
u^{\alpha}J_{\alpha}\Big).$$
This establishes (iii).

In order to prove (\ref{equal1})-(\ref{equal2}), we  observe that the geodesic
curves correspond to solutions  of the second order equations
\beq\label{geo1}
\left\{\begin{array}{l} \ddot q^i =\Gamma^i_{j,\ell}\dot q^j\dot q^\ell + 2 \Gamma^i_{j,N+\alpha}\dot q^j\dot u^\alpha + \Gamma^i_{N+\alpha,N+\beta}\dot u^\alpha\dot u^\beta\qquad i=1,\dots, N\\\,\\
 \ddot u^\gamma =\Gamma^\gamma_{j,\ell}\dot q^j\dot q^\ell + 2 \Gamma^\gamma_{j,N+\alpha}\dot q^j\dot u^\alpha + \Gamma^\gamma_{N+\alpha,N+\beta}\dot u^\alpha\dot u^\beta\qquad \gamma=1,\dots,M
 \end{array}\right.
 \eeq

By the previous construction,
for any given $w\in\R^M$ the solution of (\ref{geo1})
with initial data
\bel{gindata}
(q,u,\dot q,\dot u)(0) = (0, 0,0,w)
\eeq
satisfies
\bel{geodsol}
q(t)\equiv 0\,,\qquad\qquad u(t)= t\,w\,.
\eeq
By (\ref{geo1})
 we obtain
$$\begin{array}{l}
0~=~\ddot q^i(0)~=\Gamma^i_{N+\alpha,N+\beta}(0,0)w^\alpha w^\beta\\\,\\
0~=~\ddot u^\gamma(0)~=\Gamma^\gamma_{N+\alpha,N+\beta}(0,0)w^\alpha w^\beta
\end{array}
$$
Since these equalities hold for all initial data $w\in\R^M$ in (\ref{gindata}),
this  and property (ii) prove (\ref{Chr}).

Now, by
$$
0=\Gamma_{i,N+\alpha,N+\beta}(0,0) = ~\frac{\partial g_{i, N+\beta}}{\partial u^\gamma}\,(0,0)
+\frac{\partial g_{i, N+\gamma}}{\partial u^\beta}\,(0,0) - \frac{\partial g_{N+\beta,N+\gamma}}{\partial q^i}
\,(0,0)
$$ we obtain \ref{equal1}, since, by property (ii), one has
$$
\frac{\partial g_{i, N+\beta}}{\partial u^\gamma}\,(0,0) = - \frac{\partial g^{i, N+\beta}}{\partial u^\gamma}\,(0,0)\,,\quad  \frac{\partial g_{N+\beta,N+\gamma}}{\partial q^i}
\,(0,0) = - \frac{\partial g^{N+\beta,N+\gamma}}{\partial q^i}
\,(0,0)
$$
for all $i=1,\dots,N$ and $\alpha,\beta =1,\dots,M$.

Moreover, for every $\alpha,\beta,\gamma = 1,\dots,M$, one has
$$
 \frac{\partial g_{N+\alpha,N+\beta}}{\partial u^{\gamma}}  =
  \Gamma_{N+\alpha,N+\gamma,N+\beta} +
  \Gamma_{N+\beta, N+\gamma,N+\alpha} = 0, $$
  so that,
 in view of property  (ii),
 (\ref{equal2}) is proved as well.
 \endproof

\begin{remark} {\rm
As in
(\ref{matrixdef}), we can now define
$(e_{\alpha,\beta})_{\alpha,\beta=1,\ldots,M}$ as the inverse of the
sub-matrix $(g^{N+\alpha,N+\beta})_{\alpha,\beta=1,\ldots,M}$. Then,
since at $(q,u)=(0,0)$ we have $g_{r,s}=g^{r,s}=\delta_{r,s}$, it follows
\bel{222}
\frac{\partial e_{\alpha,\beta}}{\partial q^{i}}(0,0) =
\frac{\partial g_{N+\alpha,N+\beta}}{\partial q^{i}}(0,0)
\qquad\qquad\forall i=1,\dots,N,~~\alpha,\beta = 1,\dots,M\,.
\eeq
}
\end{remark}


  \vs

  \subsection{The orthogonal curvature of   the foliation}

  For any $(q,u)$ in the range of a  $\Lambda$-adapted chart,
 consider  the quantity
  \bel{tensor}
 \frac{\partial e_{\alpha,\beta}}{\partial q^i}\,
 dc^{\alpha}\otimes dc^\beta\otimes
  dq^i
\eeq

  \begin{lma}\label{curvinv} The function in {\rm (\ref{tensor}) }
  is intrinsically defined with respect to the foliation $\Lambda$.
  This means that if $(\tilde q,\tilde u)$ is a $\Lambda$-adapted chart then
\bel{curvtens}
  \frac{\partial \tilde e_{\alpha,\beta}}{\partial \tilde q^i} =
\frac{\partial u^\gamma}{\partial
  \tilde u^{\alpha}}  \frac{\partial u^\delta}{\partial
  \tilde u^{\beta}}  \frac{\partial q^j}{\partial
  \tilde  q^{i}}\frac{\partial e_{\gamma,\delta}}{\partial q^j}\,.
\eeq
  \end{lma}

  {\bf Proof.}  Since$(q,u)$ and $(\tilde q, \tilde u)$ are
  $\Lambda$-adapted, the coordinate transformation
  $(q,u)\mapsto (\tilde q, \tilde u)$ satisfies
  $\frac{\partial \tilde u}{\partial q} = 0$.
  Therefore,
  $$
  \tilde g^{N+\alpha, N+\beta} =  \frac{\partial\tilde u^{\alpha} }{\partial
  u^\gamma}  \frac{\partial \tilde u^{\beta} }{\partial
  u^\delta} g^{N+\gamma, N+\delta}.
  $$
By inverting the matrices on both sides of the above identity one obtains
  $$
  \tilde e_{\alpha,\beta} =  \frac{\partial u^\gamma}{\partial
  \tilde u^{\alpha}}  \frac{\partial u^\delta}{\partial
  \tilde u^{\beta}} e_{\gamma,\delta},
  $$
which implies (\ref{curvtens}),  because  $u = u(\tilde u)$
is independent of $\tilde q$.
\endproof

Although the quantity in (\ref{tensor})
is not a tensor in the strict sense of the word,  by
(\ref{curvtens}) it still transforms like a tensor w.r.t.~to changes
of $\Lambda$-adapted coordinates.
Hence, it is intrinsically defined in terms of the foliation.
By a slight abuse of language, we thus define (\ref{tensor})
as the  {\em  orthogonal curvature tensor of  the foliation} $\Lambda$.

According to Theorem \ref{mfit-th},
the foliation $\Lambda$ is $N$-fit for hyperimpulses
if and only if the
the corresponding orthogonal curvature is identically equal to zero.
We now give a geometric construction
which clarifies the meaning of
the coefficients $\partial e_{\alpha,\beta}/\partial q^i$ in
(\ref{tensor}), in the general case (see Figure 2).

Fix any point $(\q,\u)\in \Q\times\U$ and consider
any non-zero vector ${\bf V}\in \Delta^\perp_{(\q,u)}\,$.
Construct the geodesic curve that originates at $(\q,\u)$ with speed
${\bf V}$, namely
\bel{Vgeod}
s\mapsto \gamma_{\bf V}(s)\doteq Exp_{(\q,\u)} (s{\bf V})\,.
\eeq
Next, for each $s\not= 0$, consider the orthogonal space
$\Delta^\perp_{(\q_s, \u_s)}$ at the point $(\q_s,\u_s)=
\gamma_{\bf V}(s)$.
Assuming that  $s$ is sufficiently
small, a transversality argument yields the existence of a unique
vector ${\bf W}\in \Delta^\perp_{(\q_s, \u_s)}$ such that
\bel{curvedisp}
Exp_{(\q_s,\p_s)}{\bf  W}~=~(\hat \q_s,\u)~\in~ \Q\times\{\u\}\,.
\eeq
In other words, we are moving back to a point $(\hat q_s,\u)$
on the original
leaf $\Q\times \{\u\}$,  following a second geodesic curve.
In general, $\hat \q_s\not= \q$.
We claim that, setting $\sigma\doteq s^2$,  the map
$$\sigma\mapsto (\hat \q_{\sqrt\sigma}, \u)$$
defines a unique tangent vector ${\bf Z}({\bf V})\in T_{(\q,\u)}$\,.
Moreover,
the map ${\bf V} \mapsto {\bf Z}({\bf V})$ is a homogeneous
 quadratic map from $\Delta^\perp_{(\q,\u)}$ into the tangent space
 $T_{(\q,\u)}\Q\subset T_{(\q,\u)}(\Q\times\U)$.
In turn, this determines a unique symmetric bilinear mapping $B:
\Delta^\perp_{(\q,\u)}\otimes\Delta^\perp_{(\q,\u)}\mapsto T_{(\q,\u)}\Q$
such that $B({\bf V},{\bf V})={\bf Z}({\bf V})$, namely
\bel{bilin}
B({\bf V}_1,{\bf V}_2)\doteq \frac{1}{4}{\bf Z}({\bf V}_1+{\bf V}_2)
-\frac{1}{4}{\bf Z}({\bf V}_1-{\bf V}_2)\,.
\eeq

The relation between the bilinear mapping (\ref{bilin}) and the curvature tensor
(\ref{tensor})
can be best analyzed by using coordinates.
Consider
an orthonormal basis $(J_1,\dots,J_M)$ of
$\Delta^{\perp}(\q,\u)$, together with  local
$\U$-orthonormal coordinates $(q,u)$, constructed as in
Proposition \ref{orthoC}.
If ${\bf V}= w_1 J_1+\cdots + w_M J_M$, then by construction
the point $(\q_s,\u_s)$ has coordinates
$(0, sw)=(0, \ldots, 0, sw_1,\ldots, sw_M)$.
Let
$(\hat q_w(s),0)$ be the coordinates of the point $(\hat \q_s,\u)$,
constructed as in (\ref{curvedisp}).
We now have:

\begin{figure}
\centering
   \includegraphics[scale=0.50]{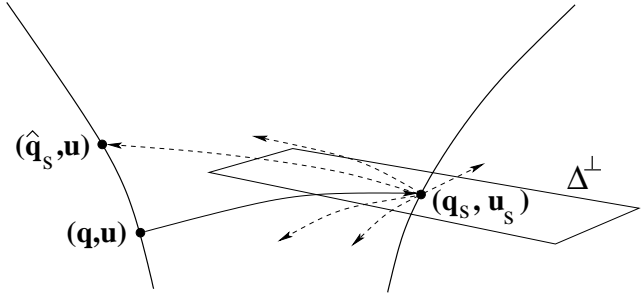}
\caption{The geodesics involved in the computation of the orthogonal
curvature of $\Lambda$.} \end{figure}

\begin{thm}\label{curvth}
The curve $s\mapsto q_w(s)\in \rr^N$ is continuous and satisfies
\bel{curvature2}
\lim_{s\to 0} \frac{\hat q_w^i(s)}{s^2} = \frac{1}{2}
\sum_{\alpha,\beta=1}^M
\frac{\partial e_{\alpha, \beta}}{\partial q^i}w^{\alpha}w^{\beta}
\qquad\qquad i=1,\ldots,N\,.
\eeq
\end{thm}

{\bf Proof}.
It is understood that the coefficients
$\partial e_{\alpha, \beta}/\partial q^i$ in (\ref{curvature2}) are
computed at $(q,u)=(0,0)$,
corresponding to the point $(\q,\u)$.
 In view of (\ref{222}), it suffices to prove that
\bel{eqform}
\lim_{s\to 0} \frac{\hat q_w^i(s)}{s^2} = \frac{1}{2}\sum_{\alpha,\beta=1}^M
\frac{\partial g_{N+\alpha, N+\beta}}{\partial q^i}w^{\alpha}w^{\beta}.
\eeq

In coordinates, the geodesic $\sigma\mapsto \gamma_{{\bf W}}(\sigma)=
 Exp_{(\q_s, \u_s)}(\sigma {\bf W})$ is given by a map
 $\sigma\mapsto (\hat q(\sigma), \hat u(\sigma))$
which, for suitable adjoint variables $p=(p_1,\ldots,p_N)$,
$\pi=(\pi_1,\ldots,\pi_M)$, satisfies the
Hamiltonian system
\beq\label{geod2}
\left\{\begin{array}{l}
\dot q^i = g^{i,j}p_j + g^{i,N+\beta}\pi_\beta\\\,
\\
\dot u^{\alpha} = g^{N+\alpha,j}p_j + g^{N+\alpha,N+\beta}\pi_\beta\\\,
\\
\dot p_i = -\frac{1}{2}\frac{\partial g^{j,k}}{\partial q^i}p_j p_k
-\frac{\partial g^{j,N+\beta}}{\partial q^i} p_j\pi_{\beta}
-\frac{1}{2}\frac{\partial g^{N+\beta,N+\gamma}}
{\partial q^i}\pi_{\beta}\pi_{\gamma}\\\,
\\
\dot \pi_{\alpha}= -\frac{1}{2}\frac{\partial g^{j,k}}{\partial u^{\alpha}}
p_jp_k  -\frac{\partial g^{j,N+\beta}}{\partial u^{\alpha}}p_j \pi_{\beta}
-\frac{1}{2}\frac{\partial g^{N+\beta,N+\gamma}}
{\partial u^{\alpha}}\pi_{\beta}\pi_{\gamma}\,.
\end{array}\right.\eeq
The conditions $\gamma_{\bf W}(0)=(\q_s, \u_s)$, $\gamma_{\bf W}(1)\in \Q\times
\{\u\}$, and the fact that ${\bf W}\in \Delta^\perp_{(\q_s,\u_s)}$ imply
\bel{bdata}
\left\{
\begin{array}{l} q^i(0)=0,\\
 u^\alpha(0)=sw_\alpha,\\
 u^\alpha(1)=0.
\end{array}\right.\qquad\qquad p_i(0)=0\,.
\eeq
For $s$ sufficiently small, the existence and uniqueness of
the solution to the two-point boundary value problem
(\ref{geod2})-(\ref{bdata}) follows from the implicit function theorem.
We now seek an expansion of this solution in powers of $s$.

Call $\bar \pi= \pi(0)$, and consider the Cauchy problem for
({\ref{geod2}), with initial data
\bel{geodindata}
\left\{\begin{array}{l}
q(0)=0\,,\\
u(0)= sw\,,
\end{array}\right.\qquad\qquad
\left\{\begin{array}{l}
p(0)=0\,,\\
\pi(0)= \bar \pi\,.
\end{array}\right.
\eeq
Using the Landau order symbols,
our computations can be simplified by observing that
\bel{orders}
\left\{\begin{array}{l}
q(\sigma)=\O(s^2)\,,\\
p(\sigma)= \O(s^2)\,,
\end{array}\right.\qquad\qquad
\left\{\begin{array}{l}
u(\sigma)=\O(s)\,,\\
\pi(\sigma)= \O(s)\,,
\end{array}\right.\qquad\qquad\forall \sigma\in [0,1]\,.
\eeq
For all $\sigma\in [0,1]$,
the solution of the Cauchy problem
(\ref{geod2}), (\ref{geodindata}) thus satisfies
\beq\label{geod3}
\left\{\begin{array}{l}\displaystyle
q^i(\sigma) =\int_0^\sigma p_i(t)\,dt +
\int_0^\sigma g^{i,N+\beta}\pi_\beta(t)\,dt
+o(s^2)\,,\\
\,
\\\displaystyle
u^{\alpha}(\sigma) = sw +\int_0^\sigma
\pi_\alpha(t)\,dt +o(s^2)\,,
\\
\,
\\
\displaystyle
p_i(\sigma) =-\frac{1}{2} \frac{\partial g^{N+\beta, N+\gamma}}{\partial q^i}
(0,0) \cdot
\int_0^\sigma \pi_\beta(t)\pi_\gamma(t)\,dt+
o(s^2)\,,
\\
\,
\\
\pi_{\alpha}(\sigma)=\bar \pi +o(s^2)\,.
\end{array}\right.\eeq

{}From the second and fourth estimates in (\ref{geod3}) we deduce
$$u^\alpha(\sigma) = sw_\alpha +\sigma\bar\pi_\alpha +o(s^2)\,.$$
Since $u^\alpha(1)=0$, this  implies
$$\bar \pi_\alpha =-sw_\alpha+o(s^2)\,.$$
Using this additional information in the third estimate, we obtain
$$
p_i(\sigma) =-\frac{1}{2} \frac{\partial g^{N+\beta, N+\gamma}}{\partial q^i}
(0,0) \cdot \sigma\, s^2 w_\beta w_\gamma+
o(s^2)\,,
$$
In turn, the first estimate now yields
$$q^i(1)=-\frac{s^2}{4} \frac{\partial g^{N+\beta, N+\gamma}}{\partial q^i}
(0,0) \, w_\beta w_\gamma - \frac{s^2}{2}
\frac{\partial g^{i, N+\beta}}{\partial u^\gamma}(0,0)\,
w_\beta w_\gamma+
o(s^2)\,.$$
Recalling the identity (\ref{equal1}), we thus obtain
$$\hat q^i(s)\doteq q^i(1)=-\frac{s^2}{2}
\frac{\partial g^{N+\beta, N+\gamma}}{\partial q^i}
(0,0) \, w_\beta w_\gamma+o(s^2)\,.$$
In view of property (ii), this establishes (\ref{eqform}).
\endproof

 \section{A variational characterization
 of input-output pairs}\label{ap3}
\setcounter{equation}{0}

Let the system $\Sigma$ be subject to
 conservative forces generated by a potential
 $U:\Q\times\U\mapsto\rr$ of class $\U^1$.
The {\em Lagrangian} of the system $(\Sigma,\g)$
subject to the potential $U$ is defined as
\bel{conslagrangian}
\ll(\q,\v,\u,{\bf w}) \doteq \T(\q,\u,\v,{\bf w}) + U(\q,\u)
\qquad \qquad ((\q,\v),(\u,{\bf w})) \in T\Q\times T\U.
\eeq
Moreover, for a given $(\u,\w)\in T\U$,  we define the {\em
$(\u,\w)$-Lagrangian } $\ll^{\u,\w}$ by setting
$$
\ll^{\u,\w}(\q,\v) \doteq \ll(\q,\v,\u,{\bf w})\qquad \qquad
(\q,\v)\in T\Q\,.$$
Let the corresponding {\em
$(\u,\w)$-Hamiltonian } be defined by
$${I\!\!H}^{\u,\p}(\q,\v) \doteq \left(\ll^{\u,\w}\right)^*(\q,\p)
 =
\sup_{v\in T_q\Q}\left\{\langle \p,\v\rangle -
\ll^{\u,\w}(\q,\v)\right\} \quad \qquad (\q,\p)\in T^*\Q\,.
$$

We now show that the admissible input-output pairs,
introduced in Definition \ref{admissible},
can be characterized by a variational principle.

\begin{thm}\label{exreg}  Assume that the Lagrangian takes the form
(\ref{conslagrangian}).
Let $u^{\sharp}:[a,b]\mapsto\U$ be a continuously differentiable control,
and let
$t\mapsto \Big(\u^{\sharp}(t)\,,~
  (\q^{\sharp}(t),\,
  \p^{\sharp}(t)\Big)$ be an {\em  admissible  input-output pair},
so that (\ref{eqham}) holds.
Then $q^{\sharp}(\cdot)$
provides a stationary point to the integral functional
\bel{intfunctional}
\I^{[a,b]}_{\u^{\sharp},\dot\u^{\sharp}}(q)\doteq \int_a^b
\ll^{\u^{\sharp}(t),\dot\u^{\sharp}(t)}(\q(t),\dot \q(t))\,dt.
\eeq
among all maps $q(\cdot)$ in the admissible set
$${\mathcal A}\doteq \Big\{\q:[a,b]\mapsto \Q\,;
~~~q(\cdot)~\hbox{is absolutely continuous},
\quad \q(a)=\q^{\sharp}(a),~~ q(b)=\q^{\sharp}(b)\Big\}\,.
$$
In this case,
the control equation
of motion has a fully Hamiltonian form,
namely
\beq\label{hamilton}
\Big(\dot\q(t),\,    \dot\p(t)\Big) =
X_{{\hh}^{\u(t),\dot\u(t)}}\Big(\q(t),\p(t)\Big) \eeq
\end{thm}
{\bf Proof.}  It is straightforward to check  that
the control equation of motion (\ref{eqham})
 coincides with the Hamiltonian system
 (\ref{hamilton}).   In turn, this is equivalent to the Euler-Lagrange
 equations for the integrand in (\ref{intfunctional}).
 Hence, a solution of (\ref{eqham}) provides a stationary
 point to the integral functional (\ref{intfunctional}).
 \endproof

\begin{remark} {\rm
Since the kinetic energy is strictly positive definite,
 the trajectory $q^{\sharp}(\cdot)$
provides a local minimizer.  More precisely, there exists $\delta>0$
such that the following holds.
For every subinterval $[\tau, \tau']\subseteq [a,b]$ such that
$\tau'-\tau\leq\delta$, one has
$$\I^{[\tau,\tau']}_{\u^{\sharp},\dot\u^{\sharp}}(q^{\sharp})
=\min_{q\in{\mathcal A}_{\tau,\tau'}}
 \I^{[\tau,\tau']}_{\u^{\sharp},\dot\u^{\sharp}}(q),$$
where the minimization is taken over the set
$${\mathcal A}_{\tau,\tau'}\doteq \Big\{\q:[\tau,\tau']
\mapsto \Q\,;~~~q(\cdot)~\hbox{is absolutely continuous},
\quad \q(\tau)=\q^{\sharp}(\tau),~~ q(\tau')=\q^{\sharp}(\tau')\Big\}\,.
$$}
\end{remark}

\part{Quadratic control systems and stabilizability}\label{analisi}
\setcounter{equation}{0}

\section{Trajectories of controlled systems with quadratic impulses}
\label{sec:4}\setcounter{equation}{0}

We now investigate  general control systems of the form:
\bel{4.1} \dot x= f(x)+\sum_{{\alpha}=1}^m g_{\alpha}(x)\,\dot u^{\alpha}
+\sum_{\alpha,\beta=1}^m
h_{{\alpha,\beta}}(x)\,\dot u^{\alpha}\dot u^{\beta}\,.
\eeq
Here the state variable $x$  and the control variable $u$  take values in
$\R^n$ and in $\R^m$, respectively.  We remark that
no a priori bounds are imposed on the derivative $\dot u$.
Our main goal is to understand under which conditions the system can be
{\it stabilized to a given point $\bar x$}.
In particular, relying on the quadratic dependence on $\dot u$ of
the right-hand side of (\ref{4.1}), in Section \ref{sec:6} we shall investigate
{\it vibrational stabilization}, achieved by means of small periodic oscillations
of the control function.
In Part \ref{applicazioni meccaniche}, these results will be applied
to the stabilization of
the mechanical systems discussed in Part \ref{meccanica}.

Throughout the following we assume that
the functions $f$, $g_{\alpha}$, and $h_{{\alpha,\beta}}=h_{\beta,\alpha}$
are at least twice
continuously differentiable.
We remark that the more general system
$$\dot x= \tilde f(t,x,u)+\sum_{{\alpha}=1}^m \tilde
g_{\alpha}(t,x,u)\,\dot u^\alpha +\sum_{\alpha,\beta=1}^m
\tilde h_{{\alpha,\beta}}(t,x,u)\,\dot u^\alpha\dot u^\beta\,,$$
where the vector fields
depend also on time and on the control $u$, can be easily
rewritten in the form (\ref{4.1}). Indeed,
it suffices to work in the extended
state space $x\in\R^{1+n+m}$, introducing the additional
state variables $x^0=t$ and $x_{n+\alpha}=u^\alpha\,$, with equations
$$\dot x^0=1\,,\qquad\qquad \dot x_{n+\alpha}=
\dot u^\alpha\qquad \alpha=1,\ldots,m\,.$$

\vskip 2em

Given the initial condition
\bel{4.2}
x(0)= \check{x}\,,
\eeq
for every smooth control function $u:[0,T]\mapsto \R^m$
 one \,\, obtains\,\, a \,\,unique
\,\,solution $t\mapsto x(t;\,u)$ of the Cauchy problem
(\ref{4.1})-(\ref{4.2}).
More generally,
since the equation (\ref{4.1}) is quadratic w.r.t.~the derivative $\dot
u$,
it is natural to consider admissible controls in a set
of absolutely continuous functions $u(\cdot)$ with derivatives
in $\L^2$.   For example, for a given $K>0$, one could allow the controls to belong to
\bel{4.3} \Big\{ u:[0,T]\mapsto\R^m\,;~~
~~~\int_0^T \big|\dot u(t)\big|^2\,dt\leq K\Big\}\,.\eeq

 The main goal of the following analysis is to
provide a characterization of the
closure of this set of trajectories,
in terms of an auxiliary differential inclusion.
Let us notice that the system (\ref{4.1})
is naturally connected with the differential inclusion
\bel{di}
\dot x \in {\mathcal F}(x),
\eeq
where,
 for every $x\in \R^n$,
\bel{dibis}
{\mathcal F}(x) \doteq  \overline{co}\left
\{ f(x) + \sum_{{\alpha}=1}^m
g_{\alpha}(x)w^{\alpha} + \sum_{\alpha,\beta=1}^m
h_{{\alpha,\beta}}(x)w^{\alpha}w^{\beta}\,;\quad
(w^1,\dots,w^m)\in \R^m\,\right\}.
\eeq
Here and in the sequel, for any given subset $A$ of a topological vector space, $\ov{co}A$ denotes the closed convex hull of
 $A$.

In addition, it will be convenient to work
also in an extended  state space,
using the variable $\hat x=\pmatrix{x^0\cr x\cr}\in\R^{1+n}$.
For a given $\hat x$, consider the set
\begin{small}
\bel{dynbdb}F(\hat x)
\doteq \co\left\{ \pmatrix{1\cr f(x)\cr}(a^0)^2+
\sum_{{\alpha}=1}^m \pmatrix{0\cr g_{\alpha}(x)\cr} a^0 a^{\alpha}
+\sum_{\alpha,\beta=1}^m
\pmatrix{0\cr h_{{\alpha,\beta}}(x)\cr}a^{\alpha}a^{\beta}\,;
~~~ a^0\in [0,1],~~\sum_{{\alpha}=0}^m
(a^{\alpha})^2 = 1\right\}.\eeq\end{small}
Notice that $F$ is a convex, compact valued multifunction on
$\R^{1+n}$, Lipschitz continuous w.r.t.~the Hausdorff metric
\cite{AC}.

For a given interval $[0,S]$, the set of trajectories of
the {\em graph differential inclusion}
\bel{4.5}\frac{d}{ds} \hat x(s)
\in F(\hat x( s))\,,\qquad\qquad \hat x(0)=\pmatrix{0\cr  x^\sharp\cr}
\eeq
is a non-empty, closed, bounded subset of ${\mathcal C}\big([0,S]\,;
~\R^{1+n}\big)$.
Consider one particular solution, say
$s\mapsto \hat x(s) =\pmatrix{x^0(s)\cr x(s)}$, defined for
$s\in [0,S]$.
Assume that  $T\doteq x^0(S)>0$. Since the map $s\mapsto x^0(s)$
is
non-decreasing, it admits a generalized inverse
\bel{4.6} s=s(t)\qquad\hbox{iff}\qquad    x^0(s)=t\,.\eeq
Indeed,
for all but countably many times $t\in [0,T]$ there exists a  unique
value of the parameter $s$ such that the identity on the right of
(\ref{4.6}) holds.
We can thus
define a corresponding trajectory
\bel{4.7}
t\mapsto x(t)=x\big(s(t)\big)\in\R^n.\eeq
This map is well defined for almost all times $t\in [0,T]$.

To establish a connection between the original control system (\ref{4.1})
and the differential inclusion (\ref{4.5}), consider first a smooth
control
function $u(\cdot)$.  As in \cite{RS}, we
define a reparametrized time variable
by setting
\bel{4.8}
s(t)\doteq \int_0^t \Big( 1+\sum_{\alpha=1}^m
(\dot u^\alpha)^2(\tau)\Big)\,
d\tau\,.\eeq
Notice that the map $t\mapsto s(t)$ is strictly increasing.
The inverse map $s\mapsto t(s)$ is uniformly
Lipschitz continuous and satisfies
$${dt\over ds}=
\left(1+\sum_{\alpha=1}^m (\dot u^\alpha)^2(t)\right)^{-1}.$$
Let now $x:[0,T]\mapsto \R^n$ be a solution of (\ref{4.1})
corresponding to the
smooth control $u:[0,T]\mapsto \R^m$.  We claim that the map
$s\mapsto \hat x(s)\doteq \pmatrix{ t(s)\cr x(t(s))\cr}$ is a solution
to the differential inclusion (\ref{4.5}).
Indeed, setting
\bel{4.9}
a^0(s)\doteq {1\over  \sqrt{1+\sum_{\beta=1}^m (\dot u^\beta)^2\big(t(s)\big)}}
\,,\qquad\qquad
a^{\alpha}(s)\doteq {\dot u^\alpha\big(t(s)\big)
\over  \sqrt{1+\sum_{\beta=1}^m (\dot u^\beta)^2\big(t(s)\big)}}\qquad
{\alpha}=1,\ldots,m\,,\eeq
one has
\bel{4.10}
\left\{
\begin{array}{l} {dt\over ds} = (a^0)^2(s)\\
\,
\\
{dx\over ds}  =
f\Big(x(s)\Big)\, (a^0)^2(s)+\sum_{{\alpha}=1}^m
g_{\alpha}\Big(x(s)\Big)\, a^0(s) a^{\alpha}(s)
+\sum_{\alpha,\beta=1}^m h_{{\alpha,\beta}}\Big(x(s)\Big)\,
a^{\alpha}(s)a^{\beta}(s)\,.\end{array}
\right.
\eeq
Hence $\hat x(\cdot) = (t(\cdot), x(\cdot)) $ verifies (\ref{4.5}),
because, by (\ref{4.9}),
$$a^0(s)\in [0,1]\,,\qquad\qquad
\sum_{{\alpha}=0}^m (a^{\alpha})^2(s)\equiv 1\,.$$
Notice that
the derivatives $\dot u^\alpha$ can now be recovered as
\bel{4.9bis}
\dot u^\alpha(t)=\frac{a^\alpha(s(t))}{a^0(s(t))}\qquad\qquad
\alpha=1,\ldots,m\,.
\eeq

The following theorem
shows that every solution of the differential inclusion
(\ref{4.5})
can be approximated by smooth solutions of the original control system
(\ref{4.1}).

\begin{thm}\label{thm4.1}
Let $\hat x=(x^0,x):[0,S]\mapsto \R^{1+n}$ be a solution
to the multivalued Cauchy problem {\rm (\ref{4.5})} such that
$x^0(S)=T>0$.  Then there exists a sequence of
smooth control functions $u_\nu:[0,T]\mapsto \R^M$
such that
the corresponding solutions
$$s\mapsto \hat x_\nu(s)=\pmatrix{t_\nu(s)\cr x_\nu(s)}$$
of the equations {\rm (\ref{4.9})-(\ref{4.10})}
converge to the map $s\mapsto \hat x(s)$ uniformly on $[0,S]$.
Moreover, defining the function
$ x(t)= x(s(t))$ as in {\rm (\ref{4.7})},
we have
\bel{4.11}\lim_{\nu\to\infty} \int_0^T \big|x(t)-x_\nu(t)\big|\,dt =0\,.
\eeq
\end{thm}

{\bf Proof.} By the assumption, the extended vector fields
$$\hat f=\pmatrix{1\cr f}\,,\qquad \hat g_{\alpha}=
\pmatrix{0\cr g_{\alpha}}\,,\qquad
\hat h_{{\alpha,\beta}}=\pmatrix{0\cr
h_{{\alpha,\beta}}}$$ are Lipschitz continuous.
Consider the set of trajectories of the control system
\bel{4.12}
{d\over ds}\, \hat x= \hat f\,\cdot (a^0)^2+\sum_{{\alpha}=1}^m
\hat g_{\alpha}\,a^0 a^{\alpha}
+\sum_{\alpha,\beta=1}^m \hat h_{ji}\, a^{\alpha}a^{\beta}\,,
\qquad\quad\hat x(0)=\pmatrix{0\cr x^\sharp}\,,
\eeq
where the controls $a=(a^0,a^1,\ldots,a^m)$ satisfy the pointwise constraints
\bel{4.13}a^0(s)\in [0,1]\,,\qquad
\sum_{{\alpha}=0}^m (a^{\alpha})^2(s)=1\qquad\qquad s\in [0,S]\,.
\eeq
In the above setting, it is well known
\cite{AC} that
the set of trajectories
$$s\mapsto \hat x(s)= (x^0,x^1,\ldots,x^n)(s)$$
of (\ref{4.12})-(\ref{4.13})
is dense on the set of solutions to the differential inclusion
(\ref{4.5}).
Hence there exists a sequence of control functions
$s\mapsto a_{\nu}(s)=\big(a^0_{\nu},\ldots, a^m_{\nu}\big)(s)$,
$\nu\geq 1$,
such that the corresponding solutions $s\mapsto \hat x_\nu(s)$
of (\ref{4.12}) converge to $\hat x(\cdot)$ uniformly for $s\in [0,S]$.
In particular, this implies the convergence of the first components:
\bel{4.14}
{x^0}_\nu(S)~=~\int_0^{S}
\big[a^0_\nu(s)\big]^2\,ds ~\to~ x^0(S) ~=~ T\,.
\eeq
We now observe that the ``input-output map''  $a(\cdot)\mapsto \hat
x(\cdot, a)$
from controls to trajectories is uniformly continuous
as a map from $\L^1\big([0,S]\,;~ \R^{1+m}\big)$ into
${\mathcal C}\big([0,S]\,;~ \R^{1+n}\big)$.
By slightly modifying the controls $a_\nu$ in $\L^1$, we can
replace the sequence $a_\nu$ by a new sequence of smooth control
functions ${\tilde a}_\nu:[0,S]\mapsto \R^{1+m}$ with the following
properties:

\bel{4.15}{{\tilde a}}_\nu^0(s)>0\qquad\qquad \hbox{for all}~~s\in [0,S]\,,
~~\nu\geq 1\,.\eeq

\bel{4.16} \int_0^{S} \big[{\tilde a}^0_\nu(s)\big]^2\,ds =T
\qquad\quad \hbox{for all}~\nu\geq 1\,,\eeq

\bel{4.17}\lim_{\nu\to\infty}~\int_0^{S}
\big| {{\tilde a}}_\nu(s) - a_\nu(s)\big|\,ds ~=~0\,.
\eeq
This implies  the uniform convergence
\bel{4.18}
\lim_{\nu\to\infty}~ \big\| \hat x(\cdot, {{\tilde a}}_\nu)-\hat x(\cdot)
\big\|_{{\mathcal C}([0,S];\,\R^{1+n})}~=~0\,.
\eeq

By (\ref{4.15}),
for each $\nu\geq 1$ the map
$$
s\mapsto x^0_\nu(s)\doteq \int_0^s \big[{\tilde a}^0_\nu(s)\big]^2\,ds
$$
is strictly increasing. Therefore it has a smooth inverse $s=s_\nu(t)$.
Recalling (\ref{4.9bis}), we now define the sequence of smooth control functions
$u_\nu:[0,T]\mapsto \R^m$ by setting
$u_\nu(t) = \big(u^1_\nu,\ldots,u^m_\nu)(t)$,
with
\bel{4.19}
u^\alpha_\nu(t)= \int_0^t\frac{{\tilde a}^\alpha_\nu(s_\nu(\tau))}
{{\tilde a}^0_\nu(s_\nu(\tau))}~ d\tau \,.\eeq
By construction, the solutions $t\mapsto x_\nu(t\,;\,u_\nu)$ of
the original system (\ref{4.1}) corresponding to the controls $u_\nu$
coincide with the trajectories
$t\mapsto (x^1_\nu,\ldots, x^n_\nu)(s_\nu(t))$,
where $\hat x_\nu= (x^0_\nu, x^1_\nu,\ldots, x^n_\nu)$ is the solution of
(\ref{4.12}) with control ${{\tilde a}}_\nu=({\tilde a}^0_\nu,\ldots,
{\tilde a}^m_\nu)$.

To prove the last statement in the theorem, define the increasing
functions
$$t(s)=\int_0^s \big[{\tilde a}^0(r)\big]^2\,dr\,,\qquad\qquad t_\nu(s)=\int_0^s
\big[{\tilde a}_\nu^0(r)\big]^2\,dr\,,
$$
and let $t\mapsto s(t)$, $t\mapsto s_\nu(t)$ be their inverses, respectively.
Notice that each $s_\nu(\cdot)$ is smooth.  Moreover,
\bel{4.20}\left|{d\over ds} t(s)\right|\leq 1\,,\qquad\qquad
\left|{d\over ds} t_\nu (s)\right|\leq 1\,,
\eeq
\bel{4.21}\lim_{\nu\to\infty}\,
\int_0^T \big| s(t)-s_\nu(t)\big|\,dt~=~\lim_{\nu\to\infty}\,
\int_0^{S} \big|t(s)-t_\nu(s)\big|\,ds=0\,.\eeq
Using (\ref{4.20}), we obtain the estimate
\begin{eqnarray}
\int_0^T\big| x(t)-x_\nu(t)\big|\,dt&=&
\int_0^T \Big| x(s(t))-x_\nu(s(t))\Big|\,dt +
\int_0^T \Big| x_\nu (s(t))-x_\nu(s_\nu(t))\Big|\,dt
\nonumber\\[-1.5ex]\label{4.22}\\[-1.5ex]
&\leq &\int_0^{S} \big| x(s)-x_\nu(s)\big|\,ds +
C\cdot \int_0^T \big|s(t)-s_\nu(t)\big|\,dt\,.
\nonumber
\end{eqnarray}

Here the constant $C$ denotes an upper bound for the
derivative w.r.t.~$s$,
for example
\bel{4.23}
C\doteq \sup_x \left\{ \big| f(x)\big|+\sum_i \big|g_{\alpha}(x)\big|+
\sum_{{\alpha,\beta}} \big|h_{{\alpha,\beta}}(x)\big|\right\}\,,
\eeq
where the supremum is taken over a compact set containing
the graphs of all functions $x_\nu(\cdot)$.
By (\ref{4.18}) and (\ref{4.21}), the right hand side of (\ref{4.22})
vanishes in the limit $\nu\to\infty$.  This completes the proof of
the theorem.
\endproof

\begin{remark} {\rm For a given time interval $[0,T]$,
we are considering controls $u(\cdot)$
in the Sobolev space $W^{1,2}$.
The corresponding solutions are absolutely continuous maps,
namely they belong to $W^{1,1}$.
Now consider a sequence of control functions $u_\nu$, whose
derivatives are uniformly bounded in ${ L}^2$.
Assume that
the corresponding reparametrized trajectories $s\mapsto (t_\nu(s), \, x_\nu(s))$,
constructed as in (\ref{4.9})-(\ref{4.10}),
converge
to a path $s\mapsto (t(s),\,x(s))$, providing a solution to (\ref{4.5}).
We wish to point out that, in general,
the projection on the state space $t\mapsto x(s(t))$
{\it may well be discontinuous}.
Notice that, on the contrary, the uniform limit of the controls $t\mapsto u_\nu(t)$
must be H\"older continuous,
because of the uniform ${\bf L}^2$ bound on the derivatives.

A completely different situation arises when all the vector fields
$h_{\alpha,\beta}$ vanish identically, so that (\ref{4.1}) reduces to
\bel{imp} \dot x= f(x)+\sum_{{\alpha}=1}^m g_{\alpha}(x)\,\dot u^\alpha
\eeq
Systems of this form have been extensively studied,  see
\cite{Su1}, \cite{miller}, \cite{B-R1}, \cite{B-R2}, or the surveys
 \cite{Rampazzo3}, \cite{B1} and the references therein.
In this case, solutions can be well defined also
for general control functions $u(\cdot)$
with bounded
variation but possibly discontinuous.  We recall that, unless the Lie brackets
$[g_\alpha,g_\beta]$  vanish identically,
 one needs to assign
a ``graph completion" of the control $u(\cdot)$ in order to determine uniquely  the trajectory.
Indeed, at each time $\tau$ where $u$ has a jump, one
should also specify a continuous path joining the left state $u(\tau-)$
with the right state $u(\tau+)$.  See \cite{B-R1} for details.

}
\end{remark}
\begin{remark}\label{linear} {\rm Assume again that system (\ref{4.1}) reduces to
(\ref{imp} ), and consider a sequence of Lipschitz controls $u_\nu$ having equi-bounded derivatives on compact sets and converging to a continuous control $\tilde u$ uniformly on bounded sets.  Moreover, consider initial states
$ x^\sharp_\nu$ converging to a point $ x^\sharp$. Then, the trajectories $x_\nu$ corresponding to the controls $u_\nu$ and the initial conditions
$x_\nu(0) =  x^\sharp_\nu$ converge to the solution $\tilde x$ of (\ref{imp} )
corresponding to the control  $\tilde u$ and initial condition
$\tilde x(0)= x^\sharp$, . In particular, if $\tilde u$ is constant, the $x_\nu$ converge to the solution of
\bel{zerog} \dot x= f(x) \qquad x(0) =  x^\sharp
\eeq
As it will be seen in Subsection \ref{nonstab}, this is quite interesting in the question of stabilizability for  mechanical systems by means of small and rapid vibrations.}
\end{remark}
\section{Stabilization }
\label{sec:5}\setcounter{equation}{0}

In this section we examine various concepts of stability
for the impulsive system (\ref{4.1}) and relate them to
the weak stability of the differential inclusion
(\ref{dynbdb})-(\ref{4.5}).

\begin{definition}\label{stabilizable} We say that the control system (\ref{4.1})
is {\em stabilizable} at the point $\bar x\in\R^n$ if, for every
$\ve>0$ there exists $\delta>0$ such that the following holds.
For every initial state $ x^\sharp$ with $| x^\sharp-\bar x|\leq\delta$
there exists a smooth control function
$t\mapsto u(t)=(u^1,\ldots,u^m)(t)$ such that the corresponding trajectory
of (\ref{4.1})-(\ref{4.2}) satisfies
\bel{5.1}
|x(t,u)-\bar x|\leq\ve\qquad\qquad\forall t\geq 0\,.
\eeq
Any such  control will be called  a {\em stabilizing control}

We say that the system (\ref{4.1}) is {\em asymptotically stabilizable}
at the point $\bar x$
if a control $u(\cdot)$ can be found such that, in addition to
(\ref{5.1}), there holds
\bel{5.2}
\lim_{t\to\infty} x(t,u)=\bar x\,.
\eeq
Any such  control will be called  an {\em asymptotically stabilizing control}.
\end{definition}

\begin{remark} {\rm Notice that the point $\bar x$ needs not to be an
equilibrium point for the vector field~$f$.}
\end{remark}

\begin{remark}{\rm We require here that the stabilizing controls be smooth.
As it will become apparent in the sequel, this is hardly a restriction.
Indeed, in all cases under consideration, if a stabilizing control
$u\in W^{1,2}$ is found, by approximation one one can construct
a smooth control $\tilde u$ which is still stabilizing.
}
\end{remark}

\begin{remark} {\rm In the above definitions we are not putting any
constraint on
the control function $u:[0,\infty[\,\mapsto\R^m$.  In principle, one may
well have
$|u(t)|\to\infty$ as $t\to\infty$.    If one wishes to stabilize the
system
(\ref{4.1}) and at the same time keep the control values within
a small neighborhood of a given value $\bar u$, it suffices to consider
the stabilization problem for an augmented system, adding the variables
$x^{n+1},\ldots, x^{n+m}$ together with the equations
$$\dot x^{n+\alpha}= \dot u^\alpha\qquad\qquad \alpha=1,\ldots,m\,.$$}
\end{remark}

Similar stability concepts can be also defined for
the differential inclusion
\bel{5.3}\dot x\in F(x)\,,\eeq
see for example \cite{Smi}.
We recall that a trajectory of (\ref{5.3})
is an absolutely continuous function $t\mapsto x(t)$ which satisfies
the differential inclusion at a.e.~time $t$.

{\bf Definition 5.2.} The point $\bar x$
is {\em weakly stable} for the differential inclusion (\ref{5.3})
if, for every $\ve>0$ there exists $\delta>0$
such that the following holds.
For every initial state ${ x^\sharp}$ with $|{ x^\sharp}-\bar x|\leq\delta$
there exists a trajectory $x(\cdot)$ of (\ref{5.3})
such that
\bel{5.4}
x(0)={ x^\sharp}\,,\qquad\qquad
|x(t)-\bar x|\leq\ve\qquad\qquad\forall t\geq 0\,.
\eeq
Moreover, $\bar x$ is {\em weakly
asymptotically stable}
if, there exists a trajectory which, in addition to
(\ref{5.4}), satisfies
\bel{5.5}\lim_{t\to\infty} x(t)=\bar x\,.
\eeq

\vskip0.5truecm

In connection with the multifunction $F$ defined at (\ref{dynbdb}), we
consider a second
multifunction $F^\pro$ obtained by projecting the sets
$F(\hat x)\subset \R^{1+n}$
into the subspace $\R^n$. More precisely, we set
\bel{5.7}
F^\pro(x)\doteq \co\left\{ f(x)\,(a^0)^2+
\sum_{{\alpha}=1}^m  g_{\alpha}(x) \,a^0 a^{\alpha}
+\sum_{\alpha,\beta=1}^m h_{{\alpha,\beta}}(x)\,a^{\alpha}a^{\beta}
~;~~~~w^0\in [0,1]\,,~~~\sum_{{\alpha}=0}^m
(w^{\alpha})^2 = 1\right\}\,.\eeq
Observe that, if the vector fields
$f, g_{\alpha}\,$, and $h_{{\alpha,\beta}}$ are Lipschitz continuous,
then the multifunction
$F^\pro$ is
Lipschitz continuous with compact, convex values.
Our first result in this section is:

\begin{thm}\label{thm5.1}
The impulsive system (\ref{4.1}) is asymptotically stabilizable
at the point $\bar x$
if and only if $\bar x$ is weakly asymptotically
stable for the {\em projected graph differential
inclusion}
\bel{5.8}\frac{d}{ds}x(s)\in F^\pro(x(s))\,.
\eeq
\end{thm}

{\bf Proof.}
Let $\bar x$ be weakly asymptotically
stable for ({\ref{5.8}). Without loss of generality, we can assume $\bar x=0$.

Given $\ve>0$, choose $\delta>0$ such that, if $|{ x^\sharp}|\leq\delta$,
then there exists a trajectory $t\mapsto x(s)$
of the differential inclusion (\ref{5.8}) such that $x(0)= x^\sharp$,
$|x(s)|\leq \ve/2$  for all $t\geq 0$ and $x(s)\to 0$ as $t\to\infty$.
Using the basic approximation property stated in Theorem \ref{thm4.1},
we will construct
a smooth control $t\mapsto u(t)= (u^1,\ldots, u^m)(t)$ such that
the corresponding trajectory $x(\cdot;u)$ of (\ref{4.1})-(\ref{4.2})
satisfies
\bel{5.9}
|x(t)|\leq\ve\qquad\forall t\geq 0\,,
\qquad\qquad\lim_{t\to\infty} x(t)= 0\,.
\eeq
Define
 the decreasing sequence of positive numbers
$\ve_k\doteq \ve\,2^{-k}$.   For each $k\geq 0$, choose $\delta_k>0$
so that, whenever $|{ x^\sharp}|\leq\delta_k$, there exists a solution to
(\ref{5.8}) with
\bel{5.10}
x(0)={ x^\sharp}\,,\quad\qquad \lim_{s\to\infty} x(s)= 0\,,\qquad\quad
|x(s)|< \frac{\ve_k}{2}\quad\forall s\geq 0\,.
\eeq
Choose a sequence of strictly
positive integers $k(1)\leq k(2)\leq \cdots$, such that
\bel{5.11}\lim_{j\to\infty} k(j)=\infty\,,\qquad\qquad
\sum_{j=1}^\infty \delta_{k(j)}=\infty \,.
\eeq
Note that the second condition in (\ref{5.11}) is
certainly satisfied if the numbers
$k(j)$ grow at a sufficiently slow rate.

Assume $|{ x^\sharp}|\leq\delta_0$.
A smooth control $u$ steering the system (\ref{4.1}) from ${ x^\sharp}$
asymptotically toward the origin will be constructed by induction on $j$.
For $j=1$, let $x:[0,s_1]\mapsto \R^n$ be a trajectory
of the differential inclusion
(\ref{5.8}) such that
$$x(0)={ x^\sharp}\,,\qquad |x(s_1)|< \frac{\delta_{k(1)}}{3}\,,\qquad
|x(s)|< {\ve_0\over 2}\quad \forall s\in[0,s_1]\,.$$
By the definition of $F^\pro$, there exists a trajectory
of the differential inclusion (\ref{4.5}) having the form
$s\mapsto\hat x(s)=(x^0(s),\, x(s))$.
Notice that, in order to apply Theorem \ref{thm4.1} and approximate $x(\cdot)$
with a smooth
solution of the control system (\ref{4.1}) we would need $x^0(s_1)>0$.
This is
not yet guaranteed by the above construction.
To take care of this problem,
we define
$s_1'\doteq s_1+\delta_{k(1)}/3C$, where $C$ provides a
local upper bound for the
magnitude of the vector field $f$, as in (\ref{4.23}).
We then prolong the trajectory $\hat x(\cdot)$
to the larger interval $[0, s_1']$, by setting
$$\frac{d}{ds}\pmatrix{x^0(s)\cr x(s)\cr} =
\pmatrix{1\cr f(x)\cr}
\qquad\qquad s\in \,]s_1,\, s_1']\,.$$
This construction achieves the inequalities
$$x^0(s_1')\geq s_1'-s_1\geq \frac{\delta_{k(1)}}{3C}\,,\qquad\qquad
|x(s_1')|< \frac{2}{3} \delta_{k(1)}\,.$$

Set $\tau^1\doteq x^0(s_1')$.
By Theorem \ref{thm4.1}, there exists a smooth control
$u:[0, \tau^1]\mapsto \R^m$
such that the corresponding solution $s\mapsto (x^0(s,u), x(s,u))$
of (\ref{4.9})-(\ref{4.10})  differs from the above trajectory
by less than $\delta_{k(1)}/3$, namely
$$|x^0(s,u)- x^0(s)|< \frac{\delta_{k(1)}}{3}\,,\qquad
|x(s;u)- x(s)|< \frac{\delta_{k(1)}}{3}\qquad\forall s\in [0,s_1']\,.$$

In particular, setting $x(t,u)\doteq x(s(t),u)$ as in (\ref{4.7}),
this implies
$$
|x(\tau_1,u)|<\delta_{k(1)}\,,\qquad\qquad
|x(t,u)|< \frac{\ve_0}{2}+\frac{\delta_{k(1)}}{3}\leq \ve_0\qquad
\forall t\in [0,\tau_1]\,.$$

The construction now proceeds by induction on $j$.  Assume
that a smooth control  $u(\cdot) $ has been constructed on the
time interval
$[0, \tau_j]$, in such a way that
\bel{5.12}
|x(\tau_j,u)|< \delta_{k(j)}\,,\qquad\qquad
|x(t,u)|< \ve_{k(j-1)}\qquad
\forall t\in [\tau_{j-1},\,\tau_j]\,.
\eeq
By assumptions, there exists a trajectory $s\mapsto x(s)$
of the differential inclusion
(\ref{5.8}) such that
\bel{5.13}
x(0)=x(\tau_j,u)\,,\qquad\qquad |x(s_j)|<\frac{\delta_{k(j+1)}}{3}\,,
\qquad\qquad |x(s)|<{\ve_{k(j)}\over 2}\quad\forall s\in [0,s_j]\,.
\eeq
This trajectory is extended to the slightly larger interval
$[0, s_j']$, with $s_j'=s_j+\delta_{k(j)}/3C$, by setting
\bel{5.14}
\frac{d}{ds}\pmatrix{x^0(s)\cr x(s)\cr}=
\pmatrix{1\cr f(x)\cr}
\qquad\qquad s\in \,]s_j,\, s_j']\,.
\eeq
Notice that, by (\ref{5.13}), (\ref{5.14}), and (\ref{4.23}),
we have
\bel{5.15}
x^0(s_j')\geq s_j'-s_j\geq \frac{\delta_{k(j)}}{3C}\,,\qquad\qquad
|x(s_j')|<\frac{2}{3} \delta_{k(j+1)}\,.
\eeq

Set $\tau_{j+1}\doteq \tau_j+x^0(s_j')$.
Using again Theorem \ref{thm4.1}, we can extend
the control $u:[0, \tau_j]\mapsto \R^M$
to a continuous, piecewise
smooth control defined on the larger interval $[0, \tau_{j+1}]$,
such that the corresponding solution $s\mapsto  x(s,u)$
of (\ref{4.1})-(\ref{4.2}) satisfies
\bel{5.16}
|x(\tau_{j+1},u)|< \delta_{k(j+1)}\,,\qquad\qquad
|x(t,u)|< \ve_{k(j)}\qquad
\forall t\in [\tau_j,\,\tau_{j+1}]\,.
\eeq
Notice that, at this stage, the control $u$ is obtained by piecing together
two smooth control functions, defined on the intervals $[0, \tau_j]$
and $[\tau_j, \tau_{j+1}]$ respectively. This makes $u$ continuous but possibly not
${\mathcal C}^1$ in a neighborhood of the point $\tau_j$.  To fix this problem,
we slightly modify the values of $u$ in a small neighborhood
of $\tau_j$, so that $u$ becomes smooth also at this point, while the
strict inequalities (\ref{5.16}) still hold.

Having completed the inductive steps for all $j\geq 1$
we observe that
$$\lim_{j\to\infty}\tau_j=\sum_j \frac{\delta_{k(j)}}{3C} = \infty$$
because of (\ref{5.11}).
As $t\to\infty$, by (\ref{5.16}) we have $x(t,u)\to 0$.
This shows that the impulsive system (\ref{4.1}) is asymptotically
stabilizable at the origin, proving one of the implications
stated in the theorem.

The converse implication is obvious, because
every solution of the system (\ref{4.1}) corresponding to a
smooth control yields a solution to the differential inclusion
(\ref{5.8}),
after a suitable time rescaling.

\begin{cor}\label{cor1}
Let a point $\bar x$ be  weakly asymptotically
stable for the differential inclusion  (\ref{di}).
Then the system (\ref{4.1}) is asymptotically stabilizable at  $\bar x$.
\end{cor}

{\bf Proof.} Since the point  $\bar x$ is weakly asymptotically stable for
{\rm (\ref{di})}, then it is asymptotically stable for  the differential inclusion
{\rm (\ref{5.8})},
 which, in turn, implies that  the impulsive system
{\rm (\ref{4.1}) } can be stabilized at  $\bar x$.

\subsection{Lyapunov functions}

There is an extensive literature, in the context of O.D.E's and of
control systems or differential inclusions, relating
the stability of an equilibrium state
to the existence of a Lyapunov function. We recall below the basic
definition,
in a form suitable for our applications.
For simplicity, we henceforth consider the case $\bar x=0\in\R^n$,
which of course is not restrictive.

\begin{definition}\label{liapdef} Let $x\mapsto G(x)\subset \R^n$
be a set valued function defined for $x\in\R^n$.
A scalar function $V$ defined on a
neighborhood
$\mathcal N$
of the origin
is a {\em weak Lyapunov function} for the differential inclusion
$$
\dot x\in G(x)
$$
if the following holds.

(i) ~$V$ is continuous on $\mathcal N$, and continuously differentiable on
$\mathcal N\setminus\{0\}$.

(ii) ~$V(0)=0$ while $V(x)>0$ for all $x\not= 0$,

(iii) For each $\delta>0$ sufficiently small, the sublevel set
$\{x\,;~~V(x)\leq\delta\}$ is compact.

(iv)~ At each $x\not= 0$ one has

\bel{5.6}
\inf_{y\in G(x)} \nabla V(x)\cdot y\leq 0\,.
\eeq

\end{definition}

The following theorem relates the stability of the impulsive control system
(\ref{4.1}) to the existence of a Lyapunov function for the differential inclusion
(\ref{di}).

\begin{thm}\label{thm5.2} Consider the multifunction ${\mathcal F}$ defined at
{\rm (\ref{dibis})}.
Assume that the differential inclusion {\rm (\ref{di})}
admits a Lyapunov function $V=V(x)$ defined on a neighborhood $\mathcal N$
of the origin.
Then the control system
{\rm (\ref{4.1}) }can be stabilized at the origin.
\end{thm}

\begin{remark}\label{ra} {\rm
We are here requiring that the function $V$ satisfies the conditions
(i)--(iii) in Definition \ref{liapdef}, and that for each
$x\not= 0$ there exists $z\in {\mathcal F}(x)$ such that
\bel{5.17'}\nabla V(x)\cdot z\leq 0\,.\eeq
Notice that the multifunction ${\mathcal F}$ in (\ref{dibis})
has unbounded values.
An equivalent condition, formulated in terms of the
bounded multifunction  $F$ in
(\ref{dynbdb}) is the following.

{\it
\noindent (iv$'$)~For every $x\in{\mathcal N}\setminus\{0\}$, there
exists $\hat y=(y_0, y)\in F(x)$ such that}
\bel{5.17}
\nabla V(x)\cdot y\leq 0\qquad\qquad y_0>0\,.\eeq
Notice that the set of conditions (i)-(iii) and (iv')
represents a slight strengthening of the notion of weak Lyapunov
function when this is applied to the projected graph differential
equation (\ref{5.8}). Yet,
let us  point out that the weak stability of
(\ref{5.8}) is not enough to guarantee the bility of
the control system (\ref{4.1}), so the condition $y_0>0$ in (\ref{5.17})
plays a crucial role.
Indeed, on $\R^2$, consider the constant vector fields $f=(1,0)$,
$h_{11}=(0,1)$, $h_{22}=(0,-1)$, $g_1=g_2=h_{12}=h_{21}=(0,0)$.
Then, choosing $a^0=0$, $a^1=a^2=1/\sqrt 2$ we see that
$(0,0,0)\in F(x)$ for every $x\in\R^2$.
Hence  condition $$
\inf_{y\in F(x)}\nabla V\cdot y \leq 0$$ is trivially satisfied
by any function $V$.
However, it is clear that in this case the system (\ref{4.1}) is not
stabilizable at the origin.}
\end{remark}

\begin{remark} {\rm
Theorem \ref{5.2} is somewhat weaker than its counterpart, Theorem \ref{thm5.1},
dealing with asymptotic stability.
Indeed, to prove that the impulsive control system (\ref{4.1})   is
stabilizable, we need to assume not only that the
differential inclusion (\ref{5.8}) is weakly stable, but also that there
exists a Lyapunov function.}
\end{remark}

{\bf Proof of Theorem \ref{thm5.2}.}
Given $\ve>0$, choose $\delta>0$ such that
$$V(x)\leq 2\delta\qquad\hbox{implies}\qquad |x|\leq\ve.$$
Let an initial state ${ x^\sharp}$ be given, with $V({ x^\sharp})\leq\delta$.

According to Remark \ref{ra}, for every $x\not= 0$ there exists $(y_0, y)\in F(x)$
such that (\ref{5.17}) holds.
We recall that the multifunction $F$ in (\ref{dynbdb}) is Lipschitz continuous,
with compact, convex values.
Since the set $\Omega\doteq
\{x\,;~~\delta\leq V(x)\leq 3\delta\}$
 is compact, by the continuity
of $\nabla V$ we can find
$\kappa>0$ such that, for every $x\in\Omega$, there exists
$\hat y=(y_0, y)\in F(x)$ with
$$\nabla V(x)\cdot y\leq 0\,,\qquad\qquad y_0\geq \kappa.$$

The control $u$ will be defined inductively on a sequence of the time intervals
$[\tau_{j-1},\,\tau_j]$, with $\tau_j\geq j\kappa$.
Set $\tau_0=0$. Consider the
differential inclusion
\bel{5.19}
\frac{d}{ds} \hat x(s)\in \left\{
\begin{array}{l}
F(x(s))\cap \{ (y_0,y)\,;~~\nabla V(x)\cdot y\leq 0\,,~~
y_0\geq \kappa\} \qquad \hbox{if}\quad
\delta< V(x)<2\delta\,, \\
F(x(s))\qquad \hbox{if}\qquad
 V(x)\leq \delta\quad \hbox{or}\quad V(x)\geq 2\delta\,,
\end{array}
\right.
\eeq
with initial data $\hat x(0)=(0, { x^\sharp})$.
The right- hand side of (\ref{5.19})
is an
upper semicontinuous multifunction, with nonempty compact convex values.
Therefore (see for example \cite{AC}), the Cauchy problem admits at
least one solution
$s\mapsto \hat x(s)= (x^0(s), x(s))$, defined for $s\in [0,1]$.
We observe that this solution satisfies
$$x^0(1)\geq \kappa\,,\qquad\qquad V(x(s))\leq \delta\qquad
\forall s\in [0,1]\,.$$
Hence, by Theorem \ref{thm4.1}
there exists a smooth control $u:[0,\tau_1]\mapsto \R^m$,
with $\tau_1=x^0(1)
\geq\kappa$,
such that the corresponding trajectory of (\ref{4.1})-\ref{4.2})
satisfies
$$V(x(t,u))< \frac{3}{2} \delta= 2\delta-2^{-1}\delta
\qquad\qquad \forall t\in [0,\tau_1]\,.$$

By induction, assume now that a smooth control $u(\cdot)$ has been
constructed on the interval $[0,\tau_j]$ with $\tau_j\geq \kappa\,j$,
and that the corresponding
trajectory $t\mapsto x(t,u)$
of the impulsive system (\ref{4.1})-(\ref{4.2}) satisfies
\bel{5.20}
V(x(t,u))\leq 2\delta-2^{-j}\delta\qquad\qquad t\in [0,\tau_j]\,.
\eeq
We then construct a solution $s\mapsto \hat x(s)=(x^0(s),\, x(s))$
of the differential inclusion
(\ref{5.19}) for $s\in [0,1]$, with initial data
$\hat x(0)=(0, \, x(\tau_j,u))$.   This function will satisfy
$$x^0(1)\geq\kappa\,,\qquad\qquad
V(x(s))<  2\delta-2^{-j}\delta \quad\forall\quad s\in [0,1]\,.$$
Using again Theorem \ref{thm4.1},
we can prolong the control $u$ to a larger
time
interval $[0,\,\tau_{j+1}]$, with $\tau_{j+1}-\tau_j=x^0(1)\geq \kappa$,
in such a way that
\bel{5.20'}
V(x(t,u))< 2\delta-2^{-j-1}\delta\qquad\qquad t\in [0,\,\tau_{j+1}]\,.
\eeq
At a first stage, this control $u$ will be piecewise smooth,
continuous but not
${\mathcal C}^1$ in a neighborhood of the point $\tau_j$.
By a local approximation, we can slightly change its values
in a small neighborhood of the point $\tau_j$,
making it smooth also at the point $\tau_j$, and preserving the
strict inequalities (\ref{5.20'}).

Since $\tau_j\geq k\,j$ for all $j\geq 1$, as $j\to\infty$
the induction procedure generates a smooth control function $u(\cdot)$, defined
for all  $t\geq 0$, whose corresponding trajectory satisfies
$V(x(t,u))<2\delta$ for all $t\geq 0$. This completes the proof of the
theorem.
\endproof

\section{A selection technique}
\label{sec:6}\setcounter{equation}{0}

In the previous section we proved two general results,
relating the stability of
the control  system (\ref{4.1}) to the weak stability of the
differential inclusion (\ref{di}).
A complete description of the sets ${\mathcal F}(x)$
in (\ref{dibis}) may often be very  difficult.
However, as shown in \cite{Smi},
to establish a stability property it
suffices to construct a suitable
family of smooth selections. We shall briefly
describe this approach.

Let a point $\bar x\in\R^n$ be given, and
assume that there exists a $\C^1$ selection
$$\gamma(x, \xi)\in \F_1(x) \doteq \ov{co}
\left\{\sum_{\alpha=1}^m g_\alpha(x) \, w^\alpha
+ \sum_{\alpha,\beta=1}^m
h_{{\alpha,\beta}}(x)\,w^{\alpha}w^{\beta}\,;\quad
(w^1,\dots,w^m)\in \R^m\,\right\} $$
depending on an additional parameter $\xi\in\R^d$, such that
\bel{ast}f(\bar x)+\gamma(\bar x, \bar\xi)=0\,.
\eeq
for some $\bar\xi\in\R^d$.
Assuming that $\gamma$ is defined on an entire neighborhood of
$(\bar x,\bar\xi)$,
consider the Jacobian matrices of partial derivatives computed at
$(\bar x,\bar\xi)$:
$$A\doteq \frac{\partial f}{\partial x}+\frac
 {\partial \gamma}{\partial x}\,,\qquad\qquad B\doteq \frac{\partial \gamma}
 {\partial \xi}\,.$$

 \begin{thm}\label{thm10.1} In the above setting,
if the linear system with constant coefficients
 \bel{lincc}
 \dot x= Ax+B\xi
 \eeq
is completely controllable, then the differential inclusion
(\ref{di})-(\ref{dibis}) is weakly asymptotically
stable at the point $\bar x$.
\end{thm}

We recall that the system (\ref{lincc}) is completely controllable if and only if
the matrices $A,B$ satisfy satisfy the algebraic relation
~Rank$\left[ B, \,AB,\,\ldots~,\, A^{n-1}B\right]=n$.
This guarantees that the system can be steered from any initial state to any final state,
within any given time interval \cite{BP, So2}.

To prove the theorem, consider the control system
\bel{conts2}
\dot x= f(x)+\gamma(x,\xi).\eeq
By a classical result in control theory, the above assumptions imply
that, for every point $x^\sharp$ sufficiently close to $\bar x$,
there exists a trajectory starting from $x^\sharp$ reaching $\bar x$
in finite time.    In particular, in view of (\ref{ast}),
 the system
(\ref{conts2}) is asymptotically stabilizable at the point $\bar x$.
Since all trajectories of (\ref{conts2}) are also trajectories of the
differential inclusion (\ref{di}), the result follows.
\endproof

\begin{remark}\label{rmk10.1}
{\rm
Toward the construction of smooth selections from the multifunction $\F$
we observe that each closed convex set $\F(x)$
can be equivalently written as
\bel{Fcones}\begin{array}{rl}
\F(x) &= ~\displaystyle f(x)+\ov{co}
\left\{\sum_{\alpha=1}^m g_\alpha(x) \, w^\alpha
+ \sum_{\alpha,\beta=1}^m
h_{{\alpha,\beta}}(x)\,w^{\alpha}w^{\beta}\,;\quad
(w^1,\dots,w^m)\in \R^m\,\right\}\\
&\qquad\qquad \displaystyle +
\ov{co}\left\{
\sum_{\alpha,\beta=1}^m
h_{{\alpha,\beta}}(x)\,w^{\alpha}w^{\beta}\,;\quad
(w^1,\dots,w^m)\in \R^m\,\right\}\\
&\doteq ~f(x)+\F_1(x)+\F_2(x)\,.\\
\end{array}
\eeq
Indeed, by definition we have $\F(x)=f(x)+\F_1(x)$.  To establish the identity
(\ref{Fcones}) it thus suffices to prove that
\bel{F12}
\F_1+\F_2\subseteq\F_1\,.\eeq
Since the set ${\mathcal F_1}(x)$ is convex and contains the origin,
for every $(w^1,\ldots, w^m)\in\R^m$
and $\ve\in [0,1]$ we have
$$y_\ve~\doteq~ \ve\left(\sum_{\alpha=1}^m g_\alpha(x)\,
\frac{w^\alpha}{\sqrt \ve}+\sum_{\alpha,\beta=1}^m h_{\alpha,\beta}(x)
\frac{w^\alpha w^\beta}{\ve}\right)~\in~\F_1\,.$$
Letting $\ve\to 0$ we find
\bel{yelim}\lim_{\ve\to 0+} \,y_\ve ~=~
\sum_{\alpha,\beta=1}^m
h_{{\alpha,\beta}}(x)\,w^{\alpha}w^{\beta}\,.
\eeq
Since ${\mathcal F}_1(x)$ is closed, it must contain
the right hand side of (\ref{yelim}).  This
proves the inclusion $\F_2\subseteq\F_1$.
Next, observing that $\F_2$ is a cone, for every $y_2\in\F_2$ and $\ve>0$ we
have $\ve^{-1}y_2\in\F_2\subseteq\F_1$.
Therefore, if $y_1\in\F_1$ we can write
$$y_1+y_2=\lim_{\ve\to 0+} (1-\ve)y_1 + \ve(\ve^{-1}y_2) \in \F_1$$
because $\F_1$ is closed and convex.  This proves (\ref{F12}).
}
\end{remark}

By Theorem \ref{thm10.1} and the above remark,
one may establish a stability result
be constructing suitable selections $\gamma(x,\xi)\in\F_2(x)$
from the cone $\F_2$.

\part{ Stabilization of mechanical systems}
\label{applicazioni meccaniche}

  In this part we address the question of how to use some
time-dependent holonomic constraints as controls in order to stabilize a
mechanical system to a given state.

\section{Stabilization with vibrating controls}
\label{mechstab} \setcounter{equation}{0}

For reader's convenience, we summarize the results in Section \ref{loc}.
Let $G= (g_{\RR,\SS})_{\RR,\SS=1,\dots,N+M}$ be  the matrix that
represents
the covariant inertial tensor in a given coordinate chart $(q,u)$.
In particular, the kinetic energy of the
whole system at a  state $(q,u)$
with velocity $(v,w)\in\R^{N+M}$
is given by
$$
\T = \frac{1}{2} g_{i,j}(q,u)v^iv^j + g_{i,N+\alpha}(q,u)v^i
w^{\alpha}
+\frac{1}{2} g_{N+\alpha,N+\beta}(q,u)w^{\alpha}w^{\beta}.
$$
Here and in the sequel, $i,j=1,\ldots,N$
while $\alpha,\beta=1,\ldots,M$.
By  $G^{-1}=(g^{r,s})_{r,s=1,\dots,N+M}$  we denote the  inverse of $G $.
Moreover, we consider the sub-matrices  $G_1\doteq
\left(g_{i,j}\right)$,
$(G^{-1})_2 \doteq \left(g^{N+\alpha,N+\beta}\right)$, and $  (G^{-1})_{12}
\doteq\left(g^{i,N+\alpha}\right)$.
Finally, we introduce
the matrices
\bel{E}
A =\left(a^{i,j}\right)\doteq (G_1)^{-1}\,,\qquad
E=\left(e_{\alpha,\beta}\right)
\doteq ((G^{-1})_2)^{-1}\,,\qquad K=\left(k^{i}_{\alpha}\right)
\doteq  (G^{-1})_{12}E\,.
\eeq
We recall that all the above  matrices depend on the variables $q,u$.
Concerning the external force,
our main assumption will be

 {\bf Hypothesis (A)}. {\it The force $F$ acting on the whole
 system does not explicitly depend on time, and is
  affine w.r.t.~the time derivative of  the control, so that }
\bel{Aforce1}F = F(q,p,u,w) =
 F^0(q,p,u)
 + F^1(q,p,u)\cdot w\,.\eeq

Taking the component along the manifold $\Q$, this implies
 $$F_\Q = F_\Q(q,p,u,w) =
 F^0_\Q(q,p, u)
 + F^1_\Q(q,p, u)\cdot w\,.$$
We can thus write the equations of motion in the form
\bel{Smech}
\left(\begin{array}{c}\dot q\\\,\\\dot p\\\,\\ \dot u\end{array}\right) =
  \left(\begin{array}{c} Ap\\\, \\-\frac{1}{2}p^\dagger \frac{\partial
A}{\partial
q}p
+ F^0_\Q \\\,\\
0\end{array}\right)
+\left(\begin{array}{c} { K}\\\,\\
-p^\dagger\frac{\partial { K}}{\partial q} +
F^1_\Q \\\,\\ 1_M\end{array}\right)\dot u\, +\,\,
\dot u^{\dag} \left(\begin{array}{c} {0}\\\,\\
\frac{1}{2} \frac{\partial E}{\partial
q} \\\,\\0\end{array}\right)\dot u \,.
\eeq

Our main goal is to find conditions which imply that the system (\ref{Smech}) is
stabilizable at a point $(\bar q,0,\bar u)$.
Two results will be described here. The first one relies  on suitable
smooth selections from the corresponding set-valued maps, as in
Theorem~\ref{thm10.1}.  The second one is based on the use of
Lyapunov functions.

 For each $q,u$, consider the cone
\bel{conegamma}\Gamma(q,u)\doteq\overline{co}\left\{
w^{\dag}  \frac{\partial E(q,u)}{\partial
q} w~ ;~~~~~w\in\R^M\right\}.
\eeq

Let $\xi\in\R^d$ be an auxiliary control variable,
ranging on a neighborhood of a
point $\bar \xi\in \R^d$, and consider a
 control system of the form
\bel{gammasyst}
\left\{\begin{array}{rl} \dot q&=Ap\,, \\
\, \\
\displaystyle \dot p&= F^0_\Q(q,p,\bar u)
+\gamma(q,p,\bar u,\xi)\,, \end{array}\right.
\eeq
where $\gamma$ is a suitable selection
 from the cone $\Gamma$.
It will be convenient to write (\ref{gammasyst}) in the more compact form
\bel{auxcontr}
(\dot q,\dot p)=\Phi(q,p,\bar u,\xi)\,,\eeq
regarding $(q,p)\in\R^{N+N}$ as state variables and
$\xi\in \R^{d}$ as control variable.
Assume that
\bel{nullpoint}
F^0_\Q(\bar q,0,\bar u)
+\gamma(\bar q,0,\bar u,\bar \xi)=0\,.\eeq
By (\ref{gammasyst}) this implies $\Phi(\bar q,0,\bar u,\bar \xi)=0\in\R^{2N}$.
To test the local controllability of (\ref{gammasyst}) at the equilibrium point
$(\bar q,0,\bar u,\bar \xi)$ we look at the linearized system
with constant coefficients
\bel{lincontr}\left(\begin{array}{c} \dot p\\  \dot q\end{array}\right)~=~
{\Lambda} \left(\begin{array}{c} p\\  q\end{array}\right)+
{\mathcal B}  \xi
\,,
\eeq
where
$${\Lambda}=\frac{\partial \Phi}{\partial (q,p)}\qquad\qquad
{\mathcal B}=\frac{\partial \Phi}{\partial \xi}
$$
with all  partial derivatives being computed at the point
$(\bar q,0,\bar u,\bar \xi)$. We can now state

\begin{thm}\label{mecstab1}
Assume that a smooth map
\bel{gammadefin}
(q,p,u,\xi)\mapsto\gamma(q,p,u,\xi)\in \Gamma(q,u)
\eeq
can be chosen in such a way that (\ref{nullpoint}) holds
and so that the linear system (\ref{lincontr}) is completely controllable.
Then the system (\ref{Smech}) is asymptotically stabilizable at the
point $(\bar q,0,\bar u)$.
\end{thm}

\noindent{\bf Proof.} According to Theorem \ref{thm10.1} and Remark \ref{rmk10.1},
it suffices to show that the control system
\bel{Smech2}
\left(\begin{array}{c}\dot q\\\,\\\dot p\\\,\\ \dot u\end{array}\right) =
  \left(\begin{array}{c} Ap \\\, \\-\frac{1}{2}p^\dagger \frac{\partial
A}{\partial
q}p
+ F^0_\Q \\\,\\
0\end{array}\right)
+\left(\begin{array}{c} { K}\\\,\\
-p^\dagger \frac{\partial { K}}{\partial q} +
F^1_\Q \\\,\\ 1_M\end{array}\right)w
+
w^{\dag} \left(\begin{array}{c} {0}\\\,\\
\frac{1}{2} \frac{\partial E}{\partial
q} \\\,\\0\end{array}\right)w+
\left(\begin{array}{c} {0}\\\,\\
\gamma(q,p,u,\xi) \\\,\\0\end{array}\right)
\eeq
is locally controllable at $(\bar q, 0,\bar u)$.  Notice that
in (\ref{Smech2}) the state variables are $q,p,u$, while $w,\xi$
are the controls.
Computing the Jacobian matrices of partial derivatives at the point
$(q,p,u;w,\xi)=(\bar q, 0,\bar u,0,\bar\xi)$, we obtain a linear
 system with constant coefficients, of the form
\bel{lincc2}
\left(\begin{array}{c}\dot q\\ \dot p \\ \dot u\end{array}\right) ~=~
\left(\begin{array}{ccc}\Lambda_{11}&0&0\\
\Lambda_{21}& \Lambda_{22}&\Lambda_{23}\\
0&0&0\end{array}\right) \left(\begin{array}{c}q\\ p\\
u\end{array}\right)+ \left(\begin{array}{cc}0&B_{12}\\
B_{21}&B_{22}\\
0&1_M\end{array}\right)
\left(\begin{array}{c}\xi\\ w \end{array}\right)~\doteq ~
\widetilde{\Lambda} \left(\begin{array}{c}q\\ p\\
u\end{array}\right)+ \Big( \widetilde {\mathcal B_1}
~~  \widetilde {\mathcal B_2}\Big)
\left(\begin{array}{c}\xi\\ w \end{array}\right)
\eeq

By assumption, the linear system (\ref{lincontr}) is completely controllable.
Therefore
\bel{span1}
\hbox{Rank}\,\left[ {\mathcal B},~{\Lambda}{\mathcal B},~\ldots~,~
{\Lambda}^{2N-1}{\mathcal B}\right] ~=~ {2N}\,.\eeq
We now observe that the matrices ${\Lambda},{ \mathcal B}$ at
(\ref{lincontr}) correspond to the
submatrices
\bel{lambdab}{\Lambda}= \left(\begin{array}{cc}\Lambda_{11}&0\\
\Lambda_{21}& \Lambda_{22} \end{array}\right)\,,\qquad\qquad {\mathcal B}=
\left(\begin{array}{c}0 \\ B_{21} \end{array}\right)\,.\eeq
Hence from (\ref{span1}) it follows
\bel{span2}
\hbox{span}\,\left[ \widetilde {\mathcal B_1},~\widetilde{\Lambda}
\widetilde{\mathcal B_1},~\ldots~,~
\widetilde{\Lambda}^{2N-1}\widetilde{\mathcal B_1}\right] ~=~
\left\{ \left(\begin{array}{c}q\\ p\\
0\end{array}\right) \,;~~~q\in\R^N,~p\in\R^N\right\}\,.\eeq
Adding to this subspace the subspace
generated by the columns of the matrix $\widetilde{\mathcal B_2}$,
we obtain the entire space $\R^{2N+M}$.   We thus conclude that the
linear system (\ref{lincc2}) is completely controllable.
In turn, this implies that the nonlinear system (\ref{Smech2}) is
asymptotically stabilizable at $(\bar q,0,\bar u)$, completing the proof.
\endproof

\vskip0.8truecm
By choosing a special kind of selection and relying of the particular structure
of (\ref{gammasyst}), we can deduce Corollary\ref{quant} below.
To state it, if $k$ is a positive integer such that $kM\geq N$ and
$W =(w_1,\dots,w_k)\in \R^{M\times k}$, let  us consider the $N\times k M$
matrix\bel{gramat}
 M(u,q,W)\doteq\left(\begin{array}{c}\frac{\partial e_{1,\beta}}
 {\partial q^1}w_1^\beta,\dots,\frac{\partial e_{M,\beta}}
 {\partial q^1}w_1^\beta,\,\,\,\dots\dots \dots,
 \frac{\partial e_{1,\beta}}{\partial q^1}w_k^\beta,\dots,
 \frac{\partial e_{M,\beta}}{\partial q^1}
 w_k^\beta\\\,\\ \cdots\cdots\\
 \, \\
\frac{\partial e_{1,\beta}}{\partial q^N}w_1^\beta,\dots,
\frac{\partial e_{M,\beta}}{\partial q^N}w_1^\beta,\,\,\,
\dots \dots\dots,\frac{\partial e_{1,\beta}}{\partial q^1}w_k^\beta,\dots,
\frac{\partial e_{M,\beta}}{\partial q^N}w_k^\beta\end{array}\right).
\eeq

\begin{cor} \label{quant}Let $k$ be a positive integer and assume that for a given state $(\bar q,\bar u)$ there exists a $k$-tuple $\bar W=(\bar w_1,\dots,\bar w_k) \in (R^M)^k$ such that
\bel{rankgramat}
\hbox{\rm Rank} \Big(M(\bar u,\bar q,\bar W)\Big) = N
\eeq
and
\bel{equilibrio}\left\{\begin{array}{l} (F_\M^0)^1 +
\sum_{\alpha,\beta=1}^M\frac{\partial e_{\alpha,\beta}}{\partial q^1}
\sum_{r=1}^k{\bar w}^{\alpha}_r{\bar w}^\beta_r = 0\\\,\qquad\qquad\qquad \cdots \,\\
(F_\M^0)^N + \sum_{\alpha,\beta=1}^M \frac{\partial e_{\alpha,\beta}}{\partial q^N}
\sum_{r=1}^2{\bar w}^{\alpha}_r{\bar w}^\beta_r = 0,\end{array}\right.
\eeq
where the involved functions are computed at $(q,p,u) = (\bar q,0,\bar u)$.
Then the system {\rm(\ref{Smech})} is asymptotically stabilizable at the
point $(\bar q,0,\bar u)$.
\end{cor}

{\bf Proof. } Let us observe that the matrices $\Lambda$ and $\mathcal B$ in
(\ref{lambdab}) have the following form:
\bel{matexpl1}
{\mathcal B} = \left(\begin{array}{c} 0_{N\times d} \\\,\\
\frac{\partial \gamma}{\partial \xi}\end{array}\right) \qquad \Lambda =
\left(\begin{array}{cc}{0_{N\times N}} & A\\\,\\
\frac{\partial (F+\gamma)}{\partial q} & \frac{\partial (F+\gamma)}{\partial p}
\end{array}\right)
\eeq
so that, in particular,
\bel{matexpcl2}
\Lambda{\mathcal B} = \left(\begin{array}{c} A\cdot
\frac{\partial \gamma}{\partial \xi} \\\,\\
\frac{\partial (F+\gamma)}{\partial p}\cdot
\frac{\partial \gamma}{\partial \xi}\end{array}\right)
\eeq

Let us set  $d=kM$,  $\xi = W = (w_1,\dots,w_k)$,
and
$$
\gamma_i(q,u,W)\doteq \frac{1}{2}\sum_{\ell=1}^k
\frac{\partial e_{\alpha,\beta}}{\partial q^i}w_\ell^{\alpha}w_\ell^{\beta}\qquad i=1,\dots,N
$$
Notice that, by $2$-homogeneity $\gamma = (\gamma^1,\dots,\gamma^N)$, is in fact a selection of the set-valued map $\Gamma$ defined in (\ref{conegamma}).
In view of Theorem \ref{mecstab1}, to prove the asymptotic stability it is sufficient find $$\bar \xi = \bar W$$ such (\ref{equilibrio}) holds and, moreover,
$$\hbox{\rm Rank}\,\, \left[{\mathcal B},\Lambda{\mathcal B}\right] (\bar q,0,\bar u,\bar W) = 2N.
$$
Since $A$ is a non-singular matrix, by (\ref{matexpcl2}) the latter condition is equivalent to
\bel{rankab}
\hbox{\rm Rank}\left(\frac{\partial \gamma}{\partial W}\right) (\bar q,0,\bar u,\bar W)\,=\, N.
\eeq
In turn, this coincides with  (\ref{rankgramat}), so the proof is concluded.
\endproof

\vskip1.0cm

We now describe a second approach, based on the construction
of a Lyapunov function.
Throughout the following we assume that the external force $F$ in
(\ref{Aforce1}) admits the representation
\bel{Aforce2}F = F(q,p,u,w) =
 -\frac{\partial U}{\partial (q,u)}
 + F^1(q,p,u)\cdot w\,.\eeq
in terms of a potential function $U=U(q,u)$.

\begin{definition}\label{effpot}
Given a $k$-tuple of vectors
${{W}}\doteq \{w_1,\dots, w_k\}\subset \R^M$,
the corresponding {\em asymptotic effective potential}
$(q,u) \mapsto  U_{ {{W}}}(q,u)$ is defined as
$$\begin{array}{rl}
&\displaystyle U_{ {{W}}}(q,u)\doteq
U(q, u)
 -\frac{1}{2}\sum_{\ell=1}^{k}
 w_\ell^{\dag} E(q, u)  w_\ell\\\,\\  &\qquad\qquad \displaystyle
\Big(=U(q, u)
 -\frac{1}{2}\sum_{\ell=1}^{k}
\sum_{\alpha,\beta=1}^M e_{\alpha,\beta}(q, u)
w_\ell^{\alpha}\ w^{\beta}_\ell\,\Big).
\end{array}$$
\end{definition}

\begin{thm}\label{mecstab2}.  Let  the external force ${ F}$
have the form (\ref{Aforce2}). For a given state
 $(\bar q,\bar u)$,  assume  that
 there exist a neighborhood ${\mathcal N} $ of  $(\bar q,\bar u)$ and a
 $k$-tuple ${{W}}\doteq \{w_1,\dots, w_k\}\subset \R^M$,
 as in Definition
 \ref{effpot} which, in addition, satisfy the following property:

$\bullet$
{ There exists a continuously differentiable map $u\mapsto \beta(u)$
defined on a neighborhood of $\bar u$  such that
 the function
 $$(q,u)\mapsto U_{ {{W}}}(q,u) + \beta(u)$$
has a strict local minimum  at $(q,u)=(\bar q,\bar u)$.}

Then
 the system (\ref{Smech}) is stabilizable at
 $(\bar q,0,\bar u)$.
\end{thm}

{\bf Proof.}  As in Section \ref{sec:6}, consider the symmetrized
differential inclusion corresponding to (\ref{Smech}), namely
\bel{symmech}
\left(\begin{array}{c}\dot q\\ \dot p\\ \dot z\end{array}\right) \in
\overline{co}
\left\{ \left(\begin{array}{c} Ap \\-\frac{1}{2}p^\dagger\frac{\partial
A}{\partial
q}p  -\frac{\partial U}{\partial q} \\
0\end{array}\right)
+\,\,
w^{\dag} \left(\begin{array}{c} {0}\\
\frac{1}{2} \frac{\partial E}{\partial
q} \\0\end{array}\right)w \,,\qquad w\in\R^M\right\}.
\eeq

To prove the theorem, it  suffices to show that
the point $(\bar q,0,\bar u)$ is a stable equilibrium
for the differential equation
\bel{odemech}
\left(\begin{array}{c}\dot q\\ \dot p\\ \dot u\end{array}\right) =
\left(\begin{array}{c}Ap\\
-\frac{1}{2} p^\dagger\frac{\partial
A}{\partial
q}p-\frac{\partial U_{ {{W}}}}{\partial q} \\
0
\end{array}\right)\,.
\eeq
Indeed, by the definition of $U_{W}$,
the right hand side of (\ref{odemech})
is a selection of the right hand-side of (\ref{symmech}).
Introducing the Hamiltonian
function
$$
H_{W} \doteq \frac{1}{2}pAp^{\dag} + U_{W} ,
$$
the equation (\ref{odemech}) can be written
in the following Hamiltonian form:
 \bel{odeham}
\Big( \dot q\,,~ \dot p\,,~ \dot u\Big)^\dagger =
 \left(\frac{\partial H_{W} }{\partial p} ~,~~
-\frac{\partial H_W }{\partial q}~,~~
0\right).
\eeq
Therefore the map \bel{lyap1}V(q,p,u) \doteq H_{W}(q,p,u)   + \beta(u)\eeq
is a Lyapunov function for (\ref{odemech}), from which it follows that
 $(\bar q,0,\bar z)$ is a stable equilibrium for (\ref{odemech}).
\endproof

\section{Examples}\label{examples}
\setcounter{equation}{0}

\begin{figure}[h]\label{fig:pend6}
\centering
   \includegraphics[scale=0.26]{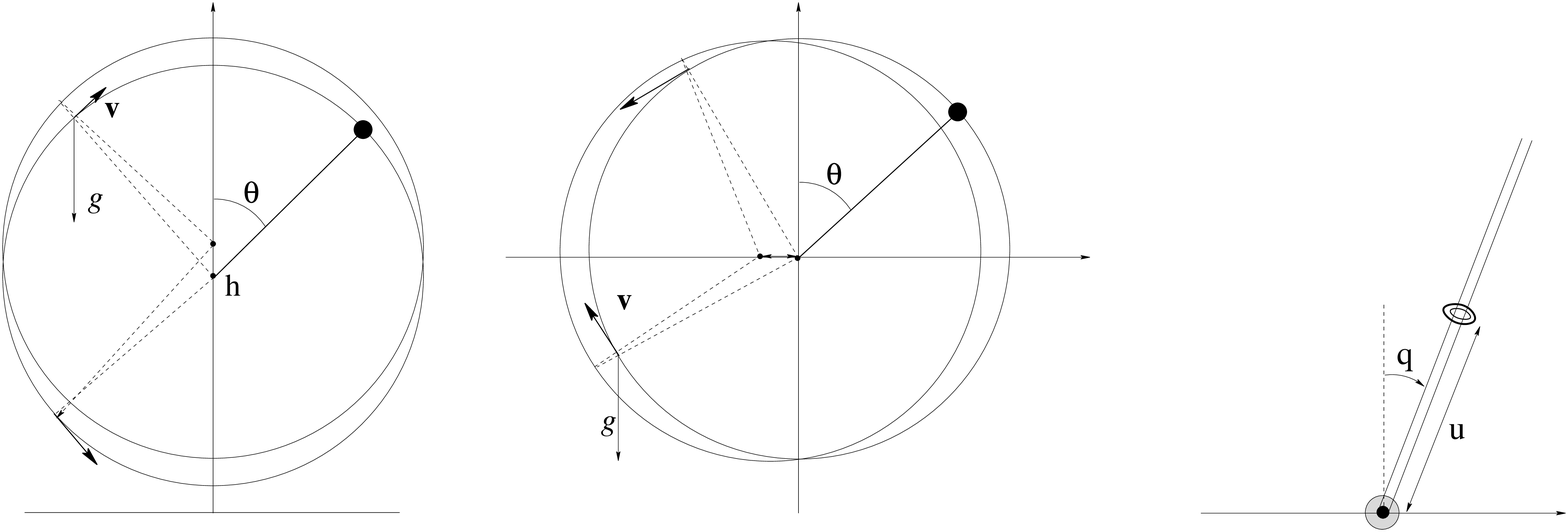}
\caption{A pendulum whose pivot oscillates vertically (on the left)
and horizontally (center). On the right: a bead sliding without friction along a rotating axis.} \end{figure}

{\bf Example 1 (pendulum with oscillating pivot).}
Let us consider a pendulum with fixed length $r=1$, whose
pivot is  moving on the vertical $y$-axis, as shown in Figure 3, left.
Its position is described by two variables: the
 clockwise
angle $\theta$ formed by the pendulum with the $y$-axis, and the height $h$
of the pivot.   We now consider $h=u(t)$ to be our control variable, while
the evolution of the other variable $\theta= q(t)$ will be
determined by the equations of motion.
We  assume that the control function $t\mapsto u(t)$
can be assigned as a function of time, ranging over a
neighborhood of the origin.
We assume that the both the pendulum and its pivot have unit mass, so that
the kinetic matrix $G$ and the matrices in
(\ref{E}) take the form
$$
G=\left(\begin{array}{lr} 1     & -\sin q \\-\sin q & 2\end{array}\right)\,
\qquad
A = (1),\qquad E = (1 + \cos^2 q),  \qquad K = (\sin q).
$$

\begin{remark}\label{pivotmass}{\rm To be consistent with the general theory
we need to put a mass on the pivot as well.  This is needed in order that
the matrix $G$ be invertible. On the other hand it is easy to show that
 the resulting control equations are independent of the mass of the pivot.
Actually this should expected, since the motion of the pivot is here considered as a
control. Of course, what is not independent of
the mass of the pivot is the constraint reaction necessary to produce a given motion of $u$.}
\end{remark}

Notice that orthogonal curvature of the constraint foliation $\Lambda$  --i.e.~the
coefficient of $(\dot u)^2$, see Section \ref{ap2}---
is different from zero, for  $\frac{dE}{dq} = -2\sin q \cos q$.

In the presence of downward
gravitational acceleration ${ g}$,
the control equations for $q$ and the corresponding momentum
$p$ is given by
\begin{equation}\label{pendulumcontr3}
\left\{\begin{array}{l} \dot q = p + (\sin q)  \dot u \\\,\\
\dot p = -\frac{\partial U}{\partial q} - p(\cos q)\dot u - (\sin q)( \cos q)
\dot u^2\, ,
\end{array}\right.
\end{equation}
 where
$U(q,u) \doteq { g}\cos q$ is the gravitational potential.

Using Theorem \ref{mecstab2}, it is easy to check that this system
is stabilizable at the upward equilibrium point $(\bar q,\bar p,\bar u)=(0,0,0)$.
Indeed, choosing $W=\{w\}$ with $w>{\it g}$, the
corresponding  effective potential
$$
U_{ W} = {\it g}\,\cos q - \frac{1}{2} (1+\cos^2 q)w^2 .
$$
has a strict local minimum at $q=0$.

To illustrate an application of Theorem \ref{mecstab1}, we now show that
the above system is asymptotically
stabilizable at every
position $(\bar q, 0,0)$ with $0<|\bar q| < \pi/2$.
To fix the ideas, assume $\bar q>0$, the other case being entirely similar.
For $\xi>0$,
the map $\gamma(q,p,\xi)=-\xi$
provides a smooth selection from the cone
$$\Gamma(q,u)~\doteq~ \ov{co}\left\{
\frac{\partial E(q,u)}{\partial  q} w^2\,;~~~
w\in\R\right\}~=~\{ -\xi\,;~~~\xi\geq 0\}.$$
The corresponding system
(\ref{gammasyst}), with $\xi $ as control variable, now takes the form
\begin{equation}\label{pendulum3}
\left\{\begin{array}{l} \dot q = p  \\
\dot p = {\it g} \, \sin q  - \xi\,.
\end{array}\right.
\end{equation}
It is easy to check that $(\bar q,\bar p,\bar\xi)= (\bar q, 0, ~{\it g}
\sin \bar q)$ is an equilibrium position and the system is locally controllable
at this point.
Indeed, the linearized control system with constant coefficients is
$$\left(\begin{array}{c}\dot  q\\ \dot  p \end{array}\right)=
\left(\begin{array}{cc} 0     & 1 \\
-g \cos \bar q &0\end{array}\right)
\left(\begin{array}{c} q\\ p \end{array}\right)+
\left(\begin{array}{c} 0 \\ -1 \end{array}\right)\xi\,.$$
By Theorem \ref{mecstab1}, the system
(\ref{pendulumcontr3}) is asymptotically stabilizable at $(\bar q,0,0)$.

By similar arguments one can show that, by means of horizontal oscillations of the pivot,
one can stabilize the system at any position of the form
$(\bar q,0,0)$, with $\frac{\pi}{2}\leq |\bar q|\leq\pi$.

\vskip5truemm

{\bf Example 2 (sliding bead).}
Consider the mechanical system  represented
in Figure~3 (right), consisting of a bead sliding without friction
along a bar, and subject to gravity.   The bar can be rotated around the
origin, in a vertical plane.
Calling $q$  the distance of the bead from the origin,
while  $u$ is the angle formed by the bar with the vertical line.
Regarding $u$ as the controlled variable,
in this case the kinetic matrix $G$ and the matrices in
(\ref{E}) take the form
$$
G=\left(\begin{array}{lr} 1     & 0 \\ 0 & q^2\end{array}\right)\,,\qquad
A=(1)\,,\qquad E=(q^2)\,,\qquad K=(0)\,.
$$
The orthogonal curvature of the constraint foliation  $\Lambda$
is not vanishing identically: indeed, one has  $\frac{dE}{dq} = 2q$.
The control equations for $q$ and the corresponding momentum
$p$ are
\bel{pendulumcontr2}
\left\{\begin{array}{l} \dot q = p\,, \\ \dot p =-{\it g}\, \cos u+ q\dot u^2\,.
\end{array}\right.
\eeq

This case is more intuitive than the previous ones.
Indeed, it is clear that a rapid oscillation of the angle $u$ generates
 a centrifugal force that can contrast the gravitational force.
More precisely, the system can be
   asymptotically stabilized at
   each $ (\bar q, \bar p, \bar  u)\in\,
    ]0,+\infty[\,\times\{0\}\times \,]-\pi/2,\pi/2[\,$.
A simple proof of this fact follows from Theorem \ref{mecstab1}.
Indeed, for $q>0$ we trivially
have $\Gamma(q,u)=\{ qw^2\,;~~w\in\R\} = \{ \xi\in\R\,;~~\xi\geq 0\}$.
It is now clear that, if $\cos \bar u>0$, then   the control system
\bel{pendulumcontr2}
\left\{\begin{array}{l} \dot q = p\,, \\ \dot p =-{\it g}\, \cos \bar u+ \xi\,,
\end{array}\right.
\eeq
admits the
 equilibrium point
$(\bar q,0, \bar \xi)$, with $\bar\xi = {\it g}\, \cos \bar u >0$.
Moreover, this system
is completely controllable around this equilibrium point, using
with controls $\xi\geq 0$.  An application of  Theorem \ref{mecstab1} yields
the asymptotic stability property.

We remark
  that the stabilizing controls cannot
  cannot be independent
   of the position $q$ and the velocity $p$.
    In particular, the approach in Theorem \ref{mecstab2},
    based on effective potential, cannot be pursued
    in this case, because  a  constant  control $w$ cannot
    stabilize the
     system
 $$
\left\{\begin{array}{l} \dot q = p\,, \\ \dot p = - {\it g} \cos u + q w^2 .
\end{array}\right.
$$

\vskip 1truecm

{\bf Example 3 (double pendulum with moving pivot).}
So far we have considered examples with  scalar controls.
We wish now to study a case where the control $u$ is two-dimensional, hence
the cone (\ref{conegamma}) is also two-dimensional.
Consider a double pendulum consisting of three point masses $P_0, P_1, P_2$,
such that the distances $|P_0P_1| ,  |P_1P_2| $ are fixed, say both equal to $1$.
Let these points be subject to the gravitational  force and constrained without friction
on a vertical  plane. Let $(u^1,u^2)$ be the cartesian coordinates of the pivot $P_0\,$,
and let $q^1,q^2$ the clockwise angles formed by $P_0P_1$
and $P_1P_2$ with the upper vertical half lines centered in $P^0$ and $P^1$,
respectively, see Figure 4.
Because of the constraints, the state of the  system $\{P_0,P_1,P_2\}$ is thus
entirely described by the four coordinates $(q^1,q^2,u^1,u^2)$.
The reduced system,
obtained by regarding the parameters  $(u^1,u^2)$ as controls and the
coordinates $(q^1,q^2)$ as
state-coordinates, is two-dimensional.
We assume that the all three points have unit mass, so that the matrix $G=(g_{rs})$
representing the kinetic energy is given by
$$G = \left(\begin{array}{cccc}
2 & \cos(q^1-q^2) & 2\cos q^1 & -2\sin q^1\\\,\\
\cos(q^1-q^2) & 1& \cos q^2 & - \sin q^2\\\,\\
2\cos q^1 &  \cos q^2 & 3 & 0\\\,\\
 -2\sin q^1 & - \sin q^2 & 0  & 3
 \end{array}\right),
 $$
Moreover, recalling (\ref{E}), we have
 $$
E = \left(\begin{array}{cc}
1 - \frac{4(\sin q^1)^2}{-3 +\cos3(q^1-q^2)}& -\frac{2\sin 2q^1}{-3 +\cos3(q^1-q^2)} \\\,\\
 -\frac{2\sin 2q^1}{-3 +\cos3(q^1-q^2)}  &1 - \frac{4(\sin q^1)^2}{-3 +\cos3(q^1-q^2)}
 \end{array}\right)\, ,
$$
$$
(F_\Q^0)^1 = 2g \sin q^1, \qquad  (F_\Q^0)^2 = g \sin q^2.
$$
Let us observe, as in Remark \ref{pivotmass}, that the matrix $E$ and the corresponding control equations are independent of the pivot's mass.
\vskip4truemm
\begin{prop}\label{bipprop}
For  every $\bar q^1\in ]0,\pi/4[$ {\rm (}resp. $\bar q^1\in ]-\pi/4,0[${\rm ) }
there exists $\delta>0$ such that for all
$\bar q^2\in ]-\delta,0[$ {\rm (}resp. $\bar q^2\in ]-\delta,0[${\rm ) }
the system is stabilizable at $(q^1,q^2,p^1,p^2,u^1,u^2) = (\bar q^1,\bar q^2, 0,0,0,0)$.\\
Moreover, the system is stabilizable at $(q^1,q^2,p^1,p^2,u^1,u^2) = (0,0, 0,0,0,0)$.
\end{prop}

\begin{figure}
\centering
   \includegraphics[scale=0.40]{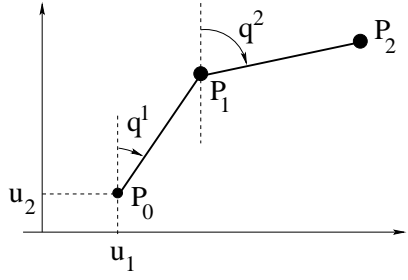}
\caption{Controlling the double pendulum by moving
the pivot at $P_0\,$.} \end{figure}

\begin{remark}{\rm For obvious reasons of translational invariance,  if we replace $(u_1,u_2)=(0,0)$ with any other value $(\bar u^1,\bar u^2)\in\R^2$} the result holds true as well.
\end{remark}
{\bf Proof of Proposition \ref{bipprop}.} Using Corollary \ref{quant} with $N=M=2$ and $k=1$, we have that the system can be stabilized to $(\bar q^1,\bar q^2,\bar u^1,\bar  u^2)$ provided there exist $\bar w\in \R^2$ such that
\bel{bipeq}
\left\{\begin{array}{l} 2g \sin \bar q^1 + \sum_{\alpha,\beta=1}^2\frac{\partial e_{\alpha,\beta}}{\partial \bar q^1} {\bar w}^{\alpha}{\bar w}^\beta = 0\\\,\\
g\sin \bar q^2 + \sum_{\alpha,\beta=1}^2 \frac{\partial e_{\alpha,\beta}}{\partial \bar q^2}{\bar w}^{\alpha}{\bar w}^\beta= 0\end{array}\right.
\eeq
and
\bel{bipdet}
\hbox{\rm det}\left(\begin{array}{cc}
 \frac{\partial e_{1,1}}{\partial q^1}\bar w^1+\frac{\partial e_{1,2}}{\partial q^1}\bar w^2 \quad
& \frac{\partial e_{2,1}}{\partial q^1}\bar w^1+\frac{\partial e_{2,2}}{\partial q^1}\bar w^2\\\,\\
 \frac{\partial e_{1,1}}{\partial q^2}\bar w^1+\frac{\partial e_{1,2}}{\partial q^2}\bar w^2 \quad
& \frac{\partial e_{2,1}}{\partial q^2}\bar w^1+\frac{\partial e_{2,2}}{\partial q^2}\bar w^2 \end{array}\right) \, \neq 0
\eeq
Notice that the latter relation can be written as
\bel{quad1}
Q_{\alpha,\beta}\bar w^\alpha \bar w^\beta \neq 0
\eeq
where the matrix $Q = \Big(Q_{\alpha,\beta}\Big)$ is defined by
\bel{formaQ}
Q\doteq \frac{\partial E}{\partial q^1}\cdot \left(\begin{array}{cc} 0~&-1\\
1~ &0\end{array}\right)\cdot  \frac{\partial E}{\partial q^2}\,.
\eeq
We recall that $E$ denotes the matrix $(e_{\alpha,\beta})$. Moreover, it  is   meant that the functions in (\ref{bipeq})-(\ref{formaQ}) are computed at $(\bar q^1,\bar q^2)$.

Let us fix $\bar q^1\in ]0,\pi/4[$. In order to establish the existence of a $\delta>0$ such that for every $\bar q^2\in ]-\delta,0[$ there is a  $\bar w$ verifying the relations (\ref{bipeq}),(\ref{bipdet}),
we need to study the intersections of the level sets of the quadratic forms $Q,\frac{\partial E}{\partial q^1}, \frac{\partial E}{\partial q^2}$.

Let us write the matrix  $\frac{\partial E}{\partial q^1}$ and $\frac{\partial E}{\partial q^2}$ explicitly:
$$\begin{array}{l}
\frac{\partial E}{\partial q^1} =  \left(\begin{array}{cc} \frac{8\sin q^1\Big(-3\cos q^1+\cos(q^1-2q^2)\Big)}{\Big(-3+\cos(2(q^1-q^2))\Big)^2}\quad &
-\frac{4\Big(-3\cos 2q^1+\cos 2q^2 \Big)}{\Big(-3+\cos(2(q^1-q^2))\Big)^2}\\\,\\
-\frac{4\Big(-3\cos 2q^1+\cos 2q^2 \Big)}{\Big(-3+\cos(2(q^1-q^2))\Big)^2}\qquad&
-\frac{8\cos q^1\Big(3\sin q^1+\sin(q^1-2q^2)\Big)}{\Big(-3+\cos(2(q^1-q^2))\Big)^2}\end{array}\right)
\\\,\\
\frac{\partial E}{\partial q^2} =  \left(\begin{array}{cc} \frac{8\sin^2 q^1 \sin(2(q^1-q^2))}{\Big(-3+\cos(2(q^1-q^2))\Big)^2}\quad &
\frac{4\sin 2q^1 \sin(2(q^1-q^2))}{\Big(-3+\cos(2(q^1-q^2))\Big)^2}\\\,\\
\frac{4\sin 2q^1 \sin(2(q^1-q^2))}{\Big(-3+\cos(2(q^1-q^2))\Big)^2}\qquad&
 \frac{8\cos^2 q^1 \sin(2(q^1-q^2))}{\Big(-3+\cos(2(q^1-q^2))\Big)^2}\end{array}\right)
\end{array}
$$
In particular, one has
$$
\hbox{\rm det}\left( \frac{\partial E}{\partial q^1}(q^1,q^2)\right) = -\frac{16}{\Big(-3+\cos(2(q^1-q^2))\Big)^2} < 0\,,\qquad \hbox{\rm det} \left(\frac{\partial E}{\partial q^2}(q^1,q^2)\right) = 0
$$
for all $q^1,q^2$.

Hence, one has:
 \begin{itemize}
 \item[(i)] The quadratic form $w\mapsto w^\dag\frac{\partial E}{\partial q^1}w$ is indefinite, so it  can be factorized by two  linear, independent, forms. Let us assume that, $\bar q^2\in]-\bar q^1,0[$, so that, in particular,  $ \frac{\partial e_{2,2}}{\partial q^1}< 0$. Hence, for suitable functions $a=a(q^1,q^2)$, $b=b(q^1,q^2)$ such that $a(q^1,q^2)\neq b(q^1,q^2)$ for all $q^1,q^2$, one has
     $$
     \frac{\partial e_{\alpha,\beta}}{\partial q^1}w^\alpha w\beta = \frac{\partial e_{2,2}}{\partial q^1}(w^2-aw^1)(w^2-bw^1).
     $$

 \item[(ii)] If  $\bar q^2\in]-\bar q^1,0[$, the quadratic form $w\mapsto w^\dag\frac{\partial E}{\partial q^2}w$ is positive semi-definite. Hence it can be factorized by  the positive scalar function $\frac{\partial e_{2,2}}{\partial q^1}$  and the square of a linear function. Moreover {\it this linear function coincides with one of the two linear  factors  of the
     quadratic form $w\mapsto w^\dag\frac{\partial E}{\partial q^1}w$}. This is a trivial consequence of the identity
     $$
     \left(\frac{\partial e_{1,2}}{\partial q^2}\frac{\partial e_{2,2}}{\partial q^1}\right)^2-2 \frac{\partial e_{2,1}}{\partial q^1}\frac{\partial e_{2,2}}{\partial q^1}\frac{\partial e_{1,2}}{\partial q^2}\frac{\partial e_{2,2}}{\partial q^1}\frac{\partial e_{2,2}}{\partial q^2} + \frac{\partial e_{1,1}}{\partial q^1}\frac{\partial e_{2,2}}{\partial q^1}\left(\frac{\partial e_{2,2}}{\partial q^2}\right)^2 = 0,
     $$
     which can be verified by direct computation. Let $(w^2-aw^1)$ be the common factor of the two quadratic forms. Hence, we obtain
     $$
     \frac{\partial e_{\alpha,\beta}}{\partial q^2}w^\alpha w^\beta = \frac{\partial e_{2,2}}{\partial q^2}(w^2-aw^1)^2.
     $$

\item[(iii)] The quadratic form  $w\mapsto w^\dag Qw$ is semi-definite and,  at each $(q^1,q^2)$, {\it it is proportional to the form $ w^\dag\frac{\partial E}{\partial q^2}w$. }More precisely, one has
    $$
    Q_{\alpha,\beta}w^{\alpha}w^\beta = \left(\frac{\partial e_{2,2}}{\partial q^1}\cdot\frac{a-b}{2} \right) \frac{\partial e_{\alpha,\beta}}{\partial q^2}w^\alpha w\beta  =
    \left( \frac{\partial e_{2,2}}{\partial q^1}\cdot\frac{\partial e_{2,2}}{\partial q^2}\cdot\frac{a-b}{2}\right) (w^2-aw^1)^2.
     $$
     This is easily  deduced by (\ref{formaQ}). Notice, in particular, the form $Q_{\alpha,\beta}w^{\alpha}w^\beta$ is never equal to the null form, since $a(q^1,q^2)\neq b(q^1,q^2)$ for all $q^1,q^2$.

     \end{itemize}

     If $S$ is a $2\times 2$ matrix and $\rho \in\R$ let us set
     $$
     \{w^\dag S w =\rho\}\doteq \{ w\in\R^2\,\,\,| \,\,\, w^\dag S w =\rho\}.
     $$
     Since $w^\dag \frac{\partial E}{\partial q^2} w$ is positive definite and  $\sin q^2<0$, there exists a real number $\eta>0$ such that
     $$
     \Big\{w^\dag \frac{\partial E}{\partial q^2} w =-\sin \bar q^2v\Big\} =
     \{w\in\R^2\,:\,\,\, (w^2-aw^1) = \eta\}\cup
     \{w\in\R^2\,:\,\,\, (w^2-aw^1) = -\eta\},
     $$
     so that, in particular,
     $$
   \Big\{w^\dag \frac{\partial E}{\partial q^2} w =-g\sin \bar q^2\Big\}\cap
   \{w\in\R^2\,:\,\,\, (w^2-aw^1) = 0\} = \emptyset.
$$
By (iii) this implies
 \bel{vuoto}
 \Big\{w^\dag \frac{\partial E}{\partial q^2} w =-g\sin \bar q^2\Big\}\cap \{w^\dag Q w =0\}
   = \emptyset.\eeq
     Moreover, by (i) the line $ \{w\in\R^2\,:\,\,\, (w^2-aw^1) = 0\}$
     is  asymptotic to the hyperbolic arc
     $$
      \Big\{w^\dag \frac{\partial E}{\partial q^1} w =-2g\sin \bar q^1\Big\} ,
      $$
which implies
\bel{soluzione}
 \Big\{w^\dag \frac{\partial E}{\partial q^1} w =-2g\sin \bar q^1\Big\}\cap
 \Big\{w^\dag \frac{\partial E}{\partial q^2} w =-g\sin \bar q^2\Big\}\neq \emptyset\,.
\eeq
Putting (\ref{vuoto}) and (\ref{soluzione}) together,
we obtain the first statement of the theorem.

On the other hand, the second statement
  will be proved by an application of
Theorem~\ref{mecstab2}. Since $U(q) = g(2\cos q^1 +\cos q^2)$
 is a potential, by letting  $W=\{(0,\eta)\}$ and $\beta(u)\doteq (u^1)^2 + (u^2)^2$,
 we have that  the effective potential
$$
U_W(q,u) \doteq U(q) + \eta^2e_{2,2}(q) + \beta(u)
$$ has a strict minimum at $(q,u) = (0,0,0,0)$ as soon as $|\eta|$ is  large enough.
In view  of Theorem~\ref{mecstab2},  this implies the that the system is stabilizable
at $(q^1,q^2,p_1,p_2,u^1,u^2) = (0,0,0,0,0,0)$.

\endproof


\section*{Acknowledgements}
The work of the first author was supported
by the N.S.F., Grant DMS-0505430.


\begin{thebibliography}{100}

\bibitem{AKN} V.~I.~Arnold, V.~V.~Kozlov, and A.~I.~Neishtadt, {\it Mathematical
Aspects of Classical and Celestial Mechanics,}
Third Edition,  Springer-Verlag, 2006.

\bibitem{AC} J.~P.~Aubin and A.~Cellina, {\it Differential
Inclusions},
Springer-Verlag, 1984.

\bibitem{B1} A.~Bressan,
Impulsive control of Lagrangian systems and locomotion in fluids,
{\it Discr. Cont. Dynam. Syst.} {\bf 20} (2008), 1-35.

\bibitem{BP} A.~Bressan and B.~Piccoli,
{\it Introduction to the Mathematical Theory of Control},
AIMS Series in Applied Mathematics, Springfield Mo. 2007.

\bibitem{B-R1} A.~Bressan and F.~Rampazzo,
On differential systems with vector-valued impulsive controls, {\it Boll.
Un. Matem. Italiana} {\bf 2-B}, (1988), 641-656.

\bibitem{B-R2} A.~Bressan and F.~Rampazzo,
Impulsive control systems with commutative vector fields, {\it
J. Optim. Theory \& Appl.} {\bf 71} (1991), 67-84.

\bibitem{B-R3} A.~Bressan and F.~Rampazzo,
On systems with quadratic impulses and their application
to Lagrangean mechanics, {\it SIAM J. Control} {\bf 31} (1993),
1205-1220.


\bibitem{AB1}
Aldo Bressan, Hyper-impulsive motions and controllizable
coordinates for Lagrangean systems {\it Atti Accad. Naz. Lincei},
Memorie, Serie VIII, Vol.~XIX (1990), 197--246.

\bibitem{AB2}
Aldo Bressan, On some control problems concerning the ski
or swing, {\it Atti Accad. Naz. Lincei, Memorie,} Serie IX, Vol.~I,
(1991), 147-196.

 \bibitem{motta}   Aldo~Bressan and  M.~Motta,  A class of mechanical systems
 with some coordinates as controls. A reduction of certain optimization problems for them. Solution methods. {\it  Atti Accad. Naz. Lincei Cl. Sci. Fis. Mat. Natur. Mem.} {\bf 9} ,Mat. Appl.  2-1 (1993), 5- 30.

\bibitem{BL} F.~Bullo and A.~D.~Lewis,
{\it Geometric Control of Mechanical Systems}, Springer-Verlag, 2004.

\bibitem{CF} Cardin,~F. and Favretti,~M. ``Hyper-Impulsive Motion on
Manifolds.'' \textit{Dynamics of continuous, discrete, and impulsive systems}
 {\bf 4}(1), 1998, pp.1-21


\bibitem{J} V.~Jurdjevic, {\it Geometric Control Theory}.
Cambridge University Press, Cambridge, 1997.

\bibitem{L1} M.~Levi, Geometry and physics of averaging with applications,
{\it Physica D} {\bf 132} (1999), 150-164.

\bibitem{L2} M.~Levi, Geometry of vibrational stabilization
and some applications,
{\it Int. J. Bifurc. Chaos} {\bf 15} (2005), 2747-2756.

\bibitem{LR} M.~Levi and Q.~Ren, Geodesics on vibrating surfaces and curvature
of the normal family, {\it Nonlinearity} {\bf 18} (2005), 2737-2743.

\bibitem{Lee} J.~M.~Lee, {\it Riemannian Manifolds. An introduction to Curvature},
Springer, 1997.


\bibitem{LS} W.~S.~Liu and H.~J.~Sussmann, Limits of highly oscillatory controls
and the approximation of general paths by admissible trajectories, in
{\it Proc. 30-th IEEE Conference on Decision and Control,}
IEEE Publications, New York, 1991, pp. 437--442.

\bibitem{Marle}
C.~Marle, G\'eom\'etrie des syst\`emes m\'ecaniques \`a liaisons actives,
in {\it Symplectic Geometry and Mathematical Physics},
260--287, P.~Donato, C.~Duval, J.~Elhadad, and G.~M.~Tuynman Eds.,
Birkh\"auser, Boston, 1991.

\bibitem{miller}B.M.~ Miller,  The generalized solutions of ordinary
differential equations in the impulse control problems,
{\it Journal of Mathematical Systems, Estimation and Control}, {\bf 4} (1994), 385--388.


\bibitem{NS}
H.~Nijmejer and A.J.~van der Schaft,
{\it Nonlinear Dynamical Control Systems},
Springer-Verlag, New York, 1990.


\bibitem{Rampazzo1} F.~Rampazzo, On Lagrangian systems
 with some coordinates as controls. \textit{Atti Accad.
 Naz. dei Lincei, Classe di Scienze Mat. Fis. Nat.} Serie 8, {\bf 82}, 1988,
 pp.685-695.



\bibitem{Rampazzo2}
F.~Rampazzo, On the Riemannian structure of a Lagrangean system and the
problem of adding time-dependent coordinates as controls
{\it European J. Mechanics
A/Solids} {\bf 10} (1991), 405-431.

\bibitem{Rampazzo3} F.~Rampazzo, Lie brackets and impulsive controls:
an unavoidable
connection,
in {\it Differential Geometry and Control,
Proc. Sympos. Pure Math.}, (AMS, Providence)(1999), 279-296.

\bibitem{RLect}  F.~Rampazzo,  {Lecture Notes on Control and Mechanics}.

\bibitem{RS} F.~Rampazzo and  C.~Sartori,
Hamilton-Jacobi-Bellman equations with fast gradient-dependence.
{\it Indiana Univ. Math. J.} {\bf 49} (2000), 1043-1078.



\bibitem{Re1} B.~L.~Reinhart,
Foliated manifolds with bundle-like metrics.  {\it Annals of Math.}
{\bf 69}  (1959), 119-132.
\v
\bibitem{Re2} B.~L.~Reinhart,
{\it Differential geometry of foliations.
The fundamental integrability problem.} Springer-Verlag, Berlin, 1983.
\v
\bibitem{Smi} G.~V.~Smirnov, {\it Introduction to the Theory of
Differential Inclusions},
American Mathematical Society,
Graduate Studies in Mathematics, vol. 41 (2002).

\bibitem{SS} J.~L.~Synge and A.~Schild, {\it Tensor Calculus},
Dover Publications, New York, 1978.

\bibitem{So}
E.~D.~Sontag,
A Lyapunov-like characterization of asymptotic
controllability, {\it SIAM J. Control \& Optim.} {\bf 21} (1983), 462-471.


\bibitem{So2}
E.D.~Sontag, {\it Mathematical Control Theory}, Springer Verlag, New
York, 1990.


\bibitem{Su1}
 H.~J.~Sussmann, On the gap between deterministic
and stochastic ordinary differential equations,
{\it Ann. Prob.} {\bf 6} (1978), 17-41.

\bibitem{Su2}
H.~J.~Sussmann,
Lie brackets, real analyticity and geometric control theory,
in {\it Mathematical Control Theory},
R.W.~Brockett, R.S.~Millmam and H.J.~Sussmann eds.,
Birkh\"{a}user, Boston Inc.~(1983), pp.~1-115.

\end{thebibliography}
\end{document}